\documentclass[]{IAC} 

\usepackage{gensymb}
\usepackage[english]{babel}
\usepackage[utf8]{inputenc}
\usepackage{newtxtext,newtxmath} 
\usepackage{enumitem} 
\let\savewidebar\widebar
\let\widebar\relax
\let\savebigtimes\bigtimes
\let\bigtimes\relax
\let\savedegree\degree
\let\degree\relax
\usepackage{mathabx}
\let\degree\savedegree
\let\bigtimes\savebigtimes
\let\widebar\savewidebar

\usepackage[dvipsnames]{xcolor}
\usepackage{tikz}
\usetikzlibrary{positioning, matrix, decorations.pathreplacing, calligraphy, shapes.geometric}
\usetikzlibrary{arrows,automata}
\usepackage{tikz-3dplot}
\usepackage[hidelinks]{hyperref}
\usepackage{cleveref}

\usepackage{algorithmic}
\usepackage{algorithm}
\usepackage{cite}
\usepackage{tabularray}
\usepackage{nicefrac}
\usepackage{siunitx}
\usepackage{fancyhdr}
\makeatletter
\newcommand\fs@nobottomruled{\def\@fs@cfont{\bfseries}\let\@fs@capt\floatc@ruled
  \def\@fs@pre{\hrule height.8pt depth0pt \kern2pt}%
  \def\@fs@post{}
  \def\@fs@mid{\kern2pt\hrule\kern2pt}%
  \let\@fs@iftopcapt\iftrue}
\makeatother

\floatstyle{nobottomruled}
\restylefloat{algorithm}

\makeatletter
\newcommand\fs@noruled{\def\@fs@cfont{\bfseries}\let\@fs@capt\floatc@ruled
  \def\@fs@post{}
  \let\@fs@iftopcapt\iftrue}
\makeatother

\addto\captionsenglish{%

}

\begin{document}

\IACpaperyear{24}
\IACpapernumber{C1.7.1}
\IACconference{75}
\IACcopyright{$75^{\textrm{th}}$}{2024}{Chiara Pozzi}{Mauro Pontani}
\IAClocation{Milan, Italy}

\title{Optimal Low-Thrust Orbit Transfers Connecting Gateway with Earth and Moon}
%


\IACauthor{Chiara Pozzi}{Department of Aerospace Engineering, Khalifa University of Science and Technology, Abu Dhabi, P.O. Box 127788, United Arab Emirates; 100064456@ku.ac.ae}
\IACauthor{Mauro Pontani}{Department of Astronautical, Electrical, and Energy Engineering, Sapienza Università di Roma, Rome, Italy; mauro.pontani@uniroma1.it}
\IACauthor{Alessandro Beolchi}{Department of Aerospace Engineering, Khalifa University of Science and Technology, Abu Dhabi, P.O. Box 127788, United Arab Emirates; 100064448@ku.ac.ae}
\IACauthor{Elena Fantino}{Department of Aerospace Engineering, Khalifa University of Science and Technology, Abu Dhabi, P.O. Box 127788, United Arab Emirates; elena.fantino@ku.ac.ae}

\abstract{Gateway will represent a primary logistic infrastructure in cislunar space. The identification of efficient orbit transfers capable of connecting Earth, Moon, and Gateway paves the way for enabling refurbishment, servicing, and utilization of this orbiting platform. This study is devoted to determining two-way minimum-time low-thrust orbit transfers that connect both Earth and Moon to Gateway. Backward time propagation is proposed as a very convenient option for trajectories that approach and rendezvous with Gateway, and a unified formulation is introduced for minimum-time orbit transfers, using either forward or backward propagation. Two-way transfers between Gateway and a specified low-altitude lunar orbit are first determined, using an indirect heuristic method, which employs the necessary conditions for optimality and a heuristic algorithm. Second, two-way orbit transfers that connect Earth and Gateway are addressed. Because these trajectories exist under the influence of two major attracting bodies, the underlying optimal control problem is formulated as a multi-arc trajectory optimization problem, involving three different representations for the spacecraft state. Multi-arc optimal control problems are associated with several corner conditions to be enforced at the junction time that separates distinct arcs. However, these conditions are shown to be solvable sequentially for the problem at hand, leveraging implicit costate transformation. This implies that the multi-arc problem has a set of unknown quantities with the same size as that of a single-arc optimal control problem, with apparent advantages of computational nature. The indirect heuristic method is also applied in this second mission scenario. Low-thrust orbit dynamics is propagated in a high-fidelity dynamical framework, with the use of planetary ephemeris and the inclusion of the simultaneous gravitational action of Sun, Earth, and Moon, along the entire transfer paths, in all cases. The numerical results unequivocally prove that the methodology developed in this research is effective for determining two-way minimum-time low-thrust orbit transfers connecting Earth, Moon, and Gateway.

}

\maketitle


\section{INTRODUCTION} 
\indent In coming years, Gateway will play the role of a vital logistic outpost for missions to the Moon, in the context of the Artemis program, as well as for deep space exploration. This orbiting platform will be assembled and operational in the next few years and placed in a synodic-resonant near rectlinear halo orbit (NRHO), in the vicinity of the exterior libration point of the Earth-Moon system. NRHO exhibits several favorable features, such as solar eclipse avoidance and easy access to the lunar South pole, to name a few \cite{artemis_plan}. NRHO traveled by Gateway is stored in an SPK-kernel compatible with the SPICE ephemeris system developed at NASA JPL \cite{lee2018sample}. \\
\indent The design of orbit transfers toward Gateway or, more generally, in cislunar space, has attracted a strong interest among the scientific community. With this regard, different options were investigated, involving either stable and unstable orbits, with the associated dynamical structures. Howell et al. \cite{Howell1} designed low-cost transfers between the Earth–Moon and Sun-Earth systems by means of transit orbits associated with the respective manifolds. Alessi et al. \cite{Alessi} proposed a two-impulse transfer strategy between Low Earth Orbits (LEO) and Lissajous orbits in the dynamical framework of the Circular Restricted 3-Body Problem (CR3BP). Pontani et al.  \cite{PONTANI2016218} introduced a polyhedral representation technique for invariant manifolds, useful to identify optimal two-impulse and low-thrust transfers from LEO to Lyapunov orbits. Further contributions focused on orbit transfers that connect a variety of Earth orbits to libration point orbits, in particular NRHO \cite{patrick2023hybrid,singh2021, he2020fireworks}. In a recent publication, Muralidharan and Howell \cite{muralidharan2023stretching} focused on leveraging stretching directions as a methodology for departure and trajectory design applications, whereas McCarty et al. \cite{MaCart} used Monotonic Basin Hopping (MBH) to identify low-thrust transfers from a Gateway-like NRHO to a distant retrograde orbit (DRO). Most recently, Sanna et al. \cite{Sanna2024} determined optimal two-impulse transfers from Gateway to a variety of Low Lunar Orbits (LLO), associated with different inclinations. Lastly, Pozzi et al. \cite{pozzi2024} found minimum-time low-thrust orbit transfers from Gateway to LLO in a high-fidelity multibody ephemeris model, using the indirect heuristic method \cite{pontani2014optimal,pontani2015minimum}. \\ 
\indent Low-thrust propulsion emerged as a viable option for orbit transfers in the last decades, due to the inherent reduced propellant consumption. However, this favorable feature is accompanied by long transfer durations, which makes the use of this technological option suitable for unmanned missions. Pioneering work on low-thrust propulsion is due to Edelbaum \cite{edelbaum1962use}, while more recent contributions were provided by Petropoulos \cite{petropoulos2003simple}, Betts \cite{betts2000very}, and Kechichian \cite{kechichian1998low}. Many works investigated optimal low-thrust trajectories in the framework of the two-body problem, with both direct and indirect optimization {methods \cite{betts2000very,pontani2014optimal,pontani2020optimal}}. Instead, limited research focused on optimal paths in multibody environments, and the great majority of contributions employed CR3BP as the dynamical framework  \cite{pierson1994three,kluever1995optimal}. In constrast, Pérez-Palau and {Epenoy \cite{perez2018fuel}} developed an indirect optimal control approach for minimum-fuel trajectories from LEO to different lunar orbits in the Sun-Earth-Moon bicircular restricted four-body problem. Herman and {Conway \cite{herman1998optimal}} employed an ephemeris model, in conjunction with a direct optimization technique, to address minimum-fuel low-thrust Earth-Moon orbit transfers, but neglected the gravitational action of the Sun. Most recently, Beolchi et al. \cite{beolchi2024} introduced a multi-arc formulation to address LEO-to-LLO low-thrust orbit transfers in a high-fidelity multibody ephemeris model, with also the inclusion of the perturbing action due to the Sun as a third body. In that work, terrestrial and lunar modified equinoctial elements (MEE) were employed, in conjunction with Cartesian coordinates as intermediate variables, because MEE were previously proven to mitigate some hypersensitivity issues encountered in the indirect numerical solution of the underlying optimal control problem \cite{pontani2020optimal}. 

\indent The work that follows relies on and considerably extends two preceding contributions by the same authors \cite{pozzi2024,beolchi2024}. In fact, it focuses on two-way low-thrust orbit transfers that connect Gateway to either a specified LEO or a prescribed LLO, in  a high-fidelity ephemeris-based dynamical framework. To address two-way transfers, both forward and backward orbit propagation is employed, and the underlying optimal control problem is posed in a unified formulation. In particular, as a major novelty with respect to previous research in the literature, backward propagation is proposed as an effective approach for rendezvous with Gateway, either from the Earth or from the Moon. Gateway-to-LLO transfers are investigated in the framework of the pertubed two-body problem (with Earth and Sun as perturbing bodies), and are formulated as single-arc optimal control problems. Instead, Gateway-to-LEO transfers involve two dominating attracting bodies (i.e., Earth and Moon), and a multi-arc optimal control problem is formulated \cite{pontani2021optimal,pontani2022optimal} as a preliminary step that precedes the use of the indirect heuristic method. In short, this research has the following objectives:
\begin{itemize} [itemsep=-0.1cm]
\item develop and use a high-fidelity ephemeris-based dynamical framework in cislunar space, to investigate orbit transfers, while employing MEE for orbit propagations;
\item introduce backward orbit propagation for transfers that approach and rendezvous with Gateway;
\item provide a unified formulation for minimum-time low-thrust transfers, using either forward or backward propagation;
\item address minimum-time two-way low-thrust transfers that connect Gateway and LLO, by formulating a single-arc optimal control problem;
\item address minimum-time two-way low-thrust transfers that connect Gateway and LEO, by formulating a multi-arc optimal control problem;
\item develop an effective indirect heuristic methodology, without any requirement for a stratified objective function (unlike the previous contribution \cite{beolchi2024});
\item solve the above mentioned minimum-time problems, in a high-fidelity dynamical framework, using the numerical technique stated at the preceding point.
\end{itemize}

\indent This work is organized as follows. Section II is concerned with orbit dynamics using MEE, as a prerequisite for describing the spacecraft motion in a high-fidelity model. Sections III and IV address minimum-time low-thrust transfers that connect either Gateway and LLO or Gateway and LEO, respectively. Both sections derive the boundary conditions associated with the underlying optimal control problem, describe the numerical solution technique, and report the numerical results. A unified formulation for minimum-time orbit transfers using either forward or backward propagation is contained in the Appendix. Finally, Section V summarizes the main accomplishments of this research.

\section{ORBIT DYNAMICS} 

\indent This study focuses on three-dimensional spacecraft trajectory optimization in the Earth-Moon system using low-thrust propulsion, in a high-fidelity dynamical framework. More specifically, orbit transfers are addressed to connect Gateway to both LEO and LLO, while adopting the following assumptions:
\begin{enumerate}[label=(\alph*),itemsep=-0.1cm]
    \item\label{a} Sun, Earth, and Moon have spherical mass distribution; 
    \item\label{b} low thrust is steerable and throttleable and is provided by a constant-specific-impulse and thrust-limited engine; 
    \item\label{c} the positions of Sun, Earth, and Moon are retrieved from a high-fidelity ephemeris {model \cite{acton1996ancillary}}. 
\end{enumerate}

\noindent Assumption \ref{a} implies that the gravitational attraction obeys the inverse square law. This is completely reasonable as the spacecraft is relatively far from both Earth and Moon for most of the transfer.



\indent Moreover, if $m$, $m_0$, $T$, and $c$ denote respectively the spacecraft mass, its initial value, the thrust magnitude, and the constant exhaust velocity of the propulsion system, then the mass ratio $m_R = \nicefrac{m}{m_{0}}$ obeys 
\begin{equation}\label{mRdot}
    \dot{m}_R = \frac{\dot{m}}{m_{0}} = -\frac{T}{m_{0} \, c} = -\frac{u_{T}}{c}
\end{equation}
\noindent where $u_T = \nicefrac{T}{m_{0}}$. As a result, the thrust acceleration magnitude is $a_T = \nicefrac{u_T}{m_R}$. 

\indent The remainder of this section describes some useful reference frames, as well as two representations for the spacecraft position and velocity, i.e., (1) Cartesian coordinates (CC) and (2) MEE. In particular, the former variables are useful for LEO-to-Gateway transfers. 

\subsection{Reference Frames}
The Earth-Centered Inertial frame (ECI) is associated with vectrix
\begin{equation}\label{ECI}
    \underline{\underline{E}} = [{\hat{c}_1}^{E} \: \: {\hat{c}_2}^{E} \: \: {\hat{c}_3}^{E}]
\end{equation}
\noindent where unit vectors ${\hat{c}_1}^{E}$ and ${\hat{c}_2}^{E}$ lie on the Earth's mean equatorial plane at epoch J2000.0. In particular, ${\hat{c}_1}^{E}$ is the vernal axis, while ${\hat{c}_3}^{E}$ is aligned with the Earth rotation axis.\\
\indent The Moon-Centered Inertial frame (MCI) is defined in relation to the ECI-frame and is associated with vectrix
\begin{equation}\label{MCI}
    \underline{\underline{M}} = [{\hat{c}_1}^{M} \: \: {\hat{c}_2}^{M} \: \: {\hat{c}_3}^{M}]
\end{equation}
\noindent where unit vectors ${\hat{c}_1}^{M}$ and ${\hat{c}_2}^{M}$ identify the lunar mean equatorial plane, which is assumed to coincide with the ecliptic plane. In particular, ${\hat{c}_1}^{M}$ is chosen to be parallel to ${\hat{c}_1}^{E}$ (i.e., the vernal axis), while ${\hat{c}_3}^{M}$ points toward the lunar rotation axis, which is orthogonal to the ecliptic plane, following the preceding assumption. The MCI-frame is related to the ECI-frame through the ecliptic obliquity angle ($\epsilon = 23.4$ degree)
\begin{equation}\label{MCIMAT}
    \underline{\underline{M}}^T = \textbf{R}_1(\epsilon)\underline{\underline{E}}^T
\end{equation}
\noindent where $\textbf{R}_j(\xi)$ denotes an elementary counterclockwise rotation about axis $j$ by a generic angle $\xi$.\\
\indent The synodic reference frame is a non-inertial frame with origin at the center of mass of the Earth-Moon system. It is associated with vectrix
\begin{equation}\label{SYNfr}
    \underline{\underline{R}}_{SYN} = [\hat{i} \: \: \hat{j} \: \: \hat{k}]
\end{equation}
\noindent where $\hat{i}$ points from Earth to Moon, $\hat{k}$ is aligned with the lunar orbit angular momentum, and $\hat{j}$ completes the right-hand sequence of unit vectors. \\
\indent Finally, the local vertical local horizontal (LVLH) frame is related to a main attracting body $B$ and rotates together with the space vehicle. It is associated with vectrix 
\begin{equation}\label{LVLHB}
    \underline{\underline{R}}_{LVLH_B} = \begin{bmatrix} {\hat{r}_B} & {\hat{\theta}_B} & {\hat{h}_B} \end{bmatrix}
\end{equation}
\noindent where $\hat{r}_B$ is aligned with the spacecraft position vector $\underrightarrow{\boldsymbol{r}}_B$ (taken from the center of mass of $B$), $\hat{h}_B$ points toward the spacecraft orbit angular momentum, whereas $\hat{\theta}_B$ completes the right-hand sequence of unit vectors. If $\Omega$, $i$, and $\theta_t$ denote respectively the instantaneous right ascension of the ascending node (RAAN), inclination, and argument of latitude (relative to body $B$), then 
\begin{equation}\label{LVLHmat}
    \underline{\underline{R}}_{LVLH}^T = \textbf{R}_3(\theta_t)\textbf{R}_1(i)\textbf{R}_3(\Omega)\underline{\underline{N}}^T
\end{equation}
\noindent where $\underline{\underline{N}}$ identifies either the ECI- or the MCI-frame. It is apparent that because either the Earth or the Moon can be regarded as the main attracting body, two distinct local vertical local horizontal reference frames $\underline{\underline{\rm R}}_{LVLH_E}$ and $\underline{\underline{\rm R}}_{LVLH_M}$ are to be introduced in this study.


\subsection{Cartesian coordinates}



\indent The spacecraft position and velocity can be defined using the respective Cartesian components  $\underline{\underline{N}}$, i.e., $(x,y,z)$ and $(v_x,v_y,v_z)$ in a convenient inertial frame. Therefore, the dynamical state of the space vehicle is identified by the vector 
\begin{equation}
    \boldsymbol{y} = \begin{bmatrix}
        x & y & z & v_x & v_y & v_z & m_R
    \end{bmatrix}^T.
\end{equation} 

\subsection{Modified equinoctial elements}

\indent In most cases, MEE are preferred to classical orbit elements (COE) for propagations, because they avoid singularities when circular or equatorial orbits are encountered (or approached). MEE are defined in terms of the classical orbit elements, i.e., semimajor axis $a$, eccentricity $e$, inclination $i$, right ascension of the ascending node (RAAN) $\Omega$, argument of periapsis $\omega$, and true anomaly $\theta_*$, as
\begin{equation}\label{COE2MEE}
    \begin{split}
        p &= a \, \left( 1 - e^2 \right) \;\;\;\;\;\;\;\;\; l = e \, \cos{(\Omega + \omega)}\\
        m &= e \, \sin{(\Omega + \omega)} \;\;\;\;\; n = \tan{\frac{i}{2}} \cos{\Omega}\\
        s &= \tan{\frac{i}{2}} \sin{\Omega} \;\;\;\;\;\;\;\;\; q = \Omega + \omega + \theta_{*}
    \end{split}
\end{equation}
\noindent where $p$ is the semilatus rectum 
and $q$ is the true longitude. These variables are nonsingular for all trajectories, with the only exception of equatorial retrograde orbits (\begingroup
\thinmuskip=0mu
\medmuskip=0mu
\thickmuskip=0mu$i = \pi$\endgroup). Letting $x_6 {\equiv} q$ and $\boldsymbol{z} {=} \begin{bmatrix}
        x_1 & x_2 & x_3 & x_4 & x_5
    \end{bmatrix}^T {\equiv}\\ {\equiv}\begin{bmatrix}
        p & l & m & n & s
    \end{bmatrix}^T$, the Gauss equations for MEE are
\begin{equation}\label{MEEdot}
    \begin{split}
        {\boldsymbol z}' &= {\rm \textbf{G}}\left({\boldsymbol z}, x_6\right){\boldsymbol a} \\
        {x}_6' &= {\sqrt \frac{\mu}{x_1^3}} \, \eta^2 + {\sqrt \frac{x_1}{\mu}} \, \frac{x_3 \, \sin{x_6} - x_5 \, \cos{x_6}}{\eta} \, a_{h}
    \end{split}
\end{equation}
\noindent where $'$ denotes the time derivative with respect to time $t$, $\mu$ represents the gravitational parameter of the main attracting body, $\eta = 1 + x_2 \cos{x_6} + x_3 \sin{x_6}$ is an auxiliary function, $\boldsymbol{a} = \begin{bmatrix} 
    a_r & a_{\theta} & a_h
\end{bmatrix}^T$ collects the non-Keplerian acceleration components in $\underline{\underline{\rm R}}_{LVLH}$, and $\rm \textbf{G}$ is a $5 \times 3$ matrix depending on $\boldsymbol z$ and $x_6$ \cite{pontanibook},
\begingroup
\thinmuskip=0mu
\medmuskip=0mu
\thickmuskip=0mu
\begin{equation}\label{Gmat}
    \resizebox{.88\hsize}{!}{${\rm \textbf{G}}\left({\boldsymbol z}, x_6\right)=\sqrt{\frac{x_1}{\mu}}\begin{bmatrix}
        0 & \cfrac{2x_1}{\eta} & 0 \\
        \sin{x_6} & \cfrac{(\eta+1)\cos{x_6}+x_2}{\eta} & -\cfrac{x_4\sin{x_6}-x_5 \cos{x_6}}{\eta}x_3 \\
        -\cos{x_6}
&  \cfrac{(\eta+1)\sin{x_6}+x_3}{\eta} &   \cfrac{x_4\sin{x_6}-x_5 \cos{x_6}}{\eta}x_2 \\
0 & 0 & \cfrac{1+x_4^2+x_5^2}{2\eta}\cos{x_6}\\
0 & 0 & \cfrac{1+x_4^2+x_5^2}{2\eta}\sin{x_6}
\end{bmatrix}$}
\end{equation}
\endgroup
In particular, the non-Keplerian-acceleration term $\boldsymbol{a}=\boldsymbol{a}_P+\boldsymbol{a}_T$,  where $\boldsymbol{a}_P$ and $\boldsymbol{a}_T$ denote two $(3\times 1)$-vectors that respectively include the components of the perturbing acceleration and the thrust acceleration in the LVLH-frame.

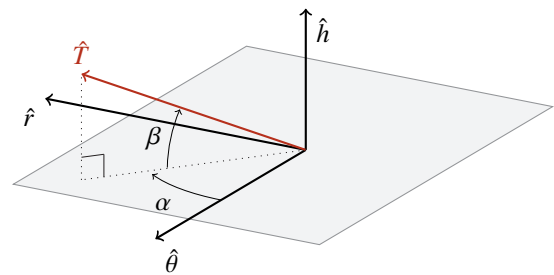
\begin{figure}[h]
    \centering
    \tdplotsetmaincoords{70}{210}
\begin{tikzpicture}[tdplot_main_coords]

\filldraw[
    draw=Gray,%
    fill=Gray!10,%
    ]          (-2,-3.1,0)
            -- (2.7,-3.1,0)
            -- (2.7,3.1,0)
            -- (-2,3.1,0)
            -- cycle;

\draw[thick,->] (0,0,0) -- (0,0,2) node[anchor=north west]{$\hat{h}$};
\draw[thick,->] (0,0,0) -- (0,4,0) node[anchor=north west]{$\hat{\theta}$};
\draw[thick,->] (0,0,0) -- (4,0,0) node[anchor=north east]{$\hat{r}$};

\pgfmathsetmacro{\ax}{2}
\pgfmathsetmacro{\ay}{2.5}
\pgfmathsetmacro{\az}{1.5}
\draw[thick,->,BrickRed] (0,0,0) -- (\ax,\ay,\az) node[anchor=south]{$\hat{T}$};


\draw[thin,dotted,Black] (0,0,0) -- (\ax,\ay,0);
\draw[thin,dotted,Black] (\ax,\ay,0) -- (\ax,\ay,\az);
\draw[very thin,Black] (\ax,\ay,0.32) -- (0.9*\ax,0.9*\ay,0.32);
\draw[very thin,Black] (0.9*\ax,0.9*\ay,0.32) -- (0.9*\ax,0.9*\ay,0);

\tdplotdrawarc[very thin,->]{(0,0,0)}{2.25}{90}{55}{anchor=north east}{$\alpha$}

\tdplotsetthetaplanecoords{53}
\tdplotdrawarc[tdplot_rotated_coords,very thin,<-]{(0,0,0)}{2}{65.5}{89}{anchor=east}{$\beta$}
 

\end{tikzpicture}
\caption{{Thrust angles $\alpha$ and $\beta$ ($0 < \alpha \leq 2\pi$, $-\nicefrac{\pi}{2} \leq \beta \leq \nicefrac{\pi}{2}$)}}
\label{Fig.ThrustAnglesAlphaBeta}
\end{figure}

\indent The dynamical state of the space vehicle is identified by the state vector 
\begin{equation}\label{stateMEE}
    \boldsymbol{x} = \begin{bmatrix}
        \boldsymbol{z}^T & x_6 & x_7
    \end{bmatrix}^T \equiv \begin{bmatrix}
        \boldsymbol{z}^T & x_6 & m_R
    \end{bmatrix}^T.
\end{equation}
\indent The thrust magnitude is closely related to $u_T$ (cf. Eq. (\ref{mRdot})), is constrained to $[0,u_T^{max}]$, and represents the first control component, i.e., $u_1:=u_T$. Moreover, in $\underline{\underline{\rm R}}_{LVLH}$, the thrust direction is identified by $\alpha$ and $\beta$ (cf. Fig. \ref{Fig.ThrustAnglesAlphaBeta}), termed in-plane and out-of-plane thrust angle, respectively. These represents two further control components. Hence, the control vector is given by 
\begin{equation}
    \boldsymbol{u} = \begin{bmatrix}
        u_1 & u_2 & u_3
    \end{bmatrix}^T \equiv \begin{bmatrix}
        u_T & \alpha & \beta
    \end{bmatrix}^T
\end{equation}
\noindent whereas $\boldsymbol{a}_T$, which collects the thrust acceleration components in $\underline{\underline{\rm R}}_{LVLH}$, is given by
\begingroup
\thinmuskip=0mu
\medmuskip=0mu
\thickmuskip=0mu
\begin{equation} \label{thrustComp}
    \boldsymbol{a}_T = \begin{bmatrix} 
    a_{T,r} & a_{T,\theta} & a_{T,h}
\end{bmatrix}^T = \frac{u_1}{x_7} \begin{bmatrix}
    \rm{s}_{u_2} \rm{c}_{u_3} & \rm{c}_{u_2} \rm{c}_{u_3} & \rm{s}_{u_3}
\end{bmatrix}^T.
\end{equation}
\endgroup
\noindent In Eq. (\ref{thrustComp}), $\rm{s}_{\theta}:=\sin\theta$ and $\rm{c}_{\theta}:=\cos\theta$ (with $\theta$ denoting a generic angle). Other than thrust, the spacecraft is affected by the gravitational pull of third bodies, i.e., Sun and Moon (when the dominating body is the Earth) or Sun and Earth (when the dominating body is the Moon). The related perturbing accelerations are projected onto the LVLH-frame and incorporated into $\boldsymbol{a}$.

\indent The state equations (\ref{mRdot}) and (\ref{MEEdot}) can be written in compact form as
\begin{equation}\label{dotx}
    {\boldsymbol{x}}' = \boldsymbol{\tilde{f}} \left( \boldsymbol{x}, \boldsymbol{u}, t \right). 
\end{equation}
These equations are accompanied by some problem-dependent boundary conditions, reported in the following sections for each specific case.

\subsection{Third-body gravitational perturbation}

\indent When a spacecraft orbits a main attracting body, the gravitational action of other massive bodies, usually referred to as third bodies, can be regarded as a perturbation acting on the spacecraft. Denoting the spacecraft with $2$, the main attracting body with $1$, the third body with $3$, and neglecting the mass of the space vehicle, the gravitational perturbation exerted on body $2$ due to body $3$ is \cite{battin1999introduction} 
\begingroup
\thinmuskip=0mu
\medmuskip=0mu
\thickmuskip=0mu
\begin{equation}\label{a3B}
    \resizebox{.88\hsize}{!}{$\underrightarrow{\boldsymbol{a}}_{3B} = \mu_3\,\left\{\frac{\underrightarrow{\boldsymbol{r}}_{13} - \underrightarrow{\boldsymbol{r}}_{12}}{\left[\left( \underrightarrow{\boldsymbol{r}}_{12} -\underrightarrow{\boldsymbol{r}}_{13} \right)\cdot\left( \underrightarrow{\boldsymbol{r}}_{12} -\underrightarrow{\boldsymbol{r}}_{13} \right)\right]^{\frac{3}{2}}} - \frac{\underrightarrow{\boldsymbol{r}}_{13}}{r_{13}^3}\right\}$}
\end{equation}
\endgroup
\noindent where {$\underrightarrow{\boldsymbol{r}}_{1j}$ is the position vector of body $j$ ($=1,2$) relative to the primary, $r_{1j}=|\underrightarrow{\boldsymbol{r}}_{1j}|$, and $\mu_3$ is the gravitational parameter of the perturbing body.} Table \ref{TABLE1} collects the gravitational parameters ${\mu}$ and the mean equatorial radii  $R$ of the relevant celestial objects.

\begin{table}[h]
\centering
\caption{Useful planetary physical parameters}
\vspace*{-1.5mm}
\begin{tblr}{Q[c,m]|Q[c,m]|Q[c,m]}
Celestial Body & ${\mu}$  $\left[\frac{\text{km}^3}{\text{s}^2}\right]$ & $R$ $\left[\text{km}\right]$ \\
\hline 
Sun & $132712440041.279$ & not used\\
Earth & $398600.436$ & $6378.136$\\ 
Moon & $4902.800$ & $1737.400$\\
\end{tblr}
\label{TABLE1}
\end{table}

\indent Because in the governing equations for MEE the non-Keplerian acceleration appears through its components in the local vertical local horizontal reference frame, both $\underrightarrow{\boldsymbol{r}}_{12}$ and $\underrightarrow{\boldsymbol{r}}_{13}$ must be projected onto $\underline{\underline{\rm R}}_{LVLH}$, 
\begin{equation}
    \boldsymbol{r}_{12}^{(LVLH)} = \begin{bmatrix}
        r_{12} \\ 0 \\ 0
    \end{bmatrix} \;\;\;\;\;\;\;\;\; \boldsymbol{r}_{13}^{(LVLH)} = \underset{\rm LVLH \leftarrow N}{\rm \textbf{R}} \, \boldsymbol{r}_{13}^{(N)}
\end{equation}
\noindent where $\boldsymbol{r}_{13}^{(N)}$ is obtained from ephemeris and N (as subscript or superscript) refers to either ECI or MCI. The preceding steps allow obtaining the components of the perturbing acceleration that appear in the column vector $\boldsymbol{a}_P$.
\section{LOW-THRUST TRANSFERS CONNECTING GATEWAY AND LLO} 

This section is focused on identifying the minimum-time transfer paths that connect Gateway and a LLO. In particular, two distinct low-thrust orbit transfers are considered: 
\begin{enumerate}[label=(\small{\alph*}),itemsep=0cm]
\vspace{-0.2cm}\item\label{GatToLLO} transfer from Gateway to LLO, and
\vspace{-0.2cm}\item\label{LLOToGat} transfer from LLO to Gateway.
\end{enumerate}

Previous research also proved that for minimum-time transfers the thrust magnitude must be set to its maximum value, i.e.,
\begin{equation}\label{u_1}
   u_1=u_T^{(max)}.
\end{equation}
This circumstance allows integrating $x_7$ in closed form. As a result, this variable can be eliminated from the state vector.

\subsection{Formulation of the problem}
\indent In case \ref{GatToLLO}, forward propagation is used, by introducing the auxiliary forward time variable $\tau_F:=t-t_0$. Both the initial time $t_0$ and the final time $t_f$ are unspecified, and $\tau_F$, which plays the role of independent variable for the dynamical system at hand, is constrained to $[0,t_f-t_0]$; $(t_f-t_0)$ represents the time of flight, also denoted with $\tau_{fin}$ hence forward. Time $t_0$ identifies the initial state (cf. Eq. (\ref{stateMEE})) at departure, which is collected in
\begin{equation}\label{Psi_0_a}
    \boldsymbol{\Psi}_i^{(a)}=  \begin{bmatrix}
       x_{1,i}-p(t_0) \\
       x_{2,i}-l(t_0)\\
       x_{3,i}-m(t_0)\\
       x_{4,i}-n(t_0)\\
       x_{5,i}-s(t_0)\\
       x_{6,i}-q(t_0)        
    \end{bmatrix} = \boldsymbol{0}.
\end{equation}
The final state corresponds to orbit injection into a circular LLO, with prescribed radius and inclination. This leads to defining the final conditions, in the form
\begin{equation}\label{Psi_f_a}
    \boldsymbol{\Psi}_f^{(a)}=  \begin{bmatrix}
       x_{1,f}-p_d \\ 
       x_{2,f}^2+x_{3,f}^2-e_d^2\\ 
       x_{4,f}^2+x_{5,f}^2-\tan^2{\cfrac{i_d}{2}}
    \end{bmatrix} = \boldsymbol{0}
\end{equation}
where $p_d=R_{LLO}$ denotes the desired semilatus rectum, $e_d$ $(=0)$ the final eccentricity, and $i_d$ the desired inclination. For case \ref{GatToLLO}, because forward propagation is employed, the state equations (\ref{dotx}) remain unaltered when rewritten with $\tau_F$ as the independent variable, i.e.,
\begin{equation}\label{stateEqforw}
     \dot{\boldsymbol{x}} =\boldsymbol{f} \left( \boldsymbol{x}, \boldsymbol{u}, \boldsymbol{p}, \tau_F \right) = \boldsymbol{\tilde{f}} \left( \boldsymbol{x}, \boldsymbol{u}, \tau_F+t_0 \right). 
\end{equation}
where the dot denotes the derivative with respect to $\tau_F$, whereas $\boldsymbol{p}$ represents a vector of unknown time-independent parameters, i.e., the two epochs $t_0$ and $t_f$ in this research. Moreover, the thrust acceleration magnitude $a_T$ is given by
\begin{equation}\label{aTforw}
     a_T=\frac{u_T^{(max)}}{m_R}=\frac{u_T^{(max)}}{1-\frac{u_T^{(max)}}{c}\tau_F}
\end{equation}

\indent In case \ref{LLOToGat}, the spacecraft travels from a circular LLO with specified radius $R_{LLO}$ and inclination $i_d$ toward Gateway. Because the final dynamical conditions, i.e., rendezvous with Gateway, are more stringent than the initial conditions, it is convenient to investigate this transfer using backward time propagation. To this end, the auxiliary backward time variable $\tau_B:=t_f-t$ is introduced. Both the initial (actual) time $t_0$ and the final actual time $t_f$ are unspecified, therefore $\tau_B$, which plays the role of independent variable for the dynamical system at hand, is monotonically increasing and constrained to $[0,t_f-t_0]$. With respect to $\tau_B$, the actual time $t_f$ identifies the initial state (cf. Eq. (\ref{stateMEE})), which is collected in
\begin{equation}\label{Psi_0_b}
    \boldsymbol{\Psi}_i^{(b)}=  \begin{bmatrix}
       x_{1,i}-p(t_f) \\
       x_{2,i}-l(t_f)\\
       x_{3,i}-m(t_f)\\
       x_{4,i}-n(t_f)\\
       x_{5,i}-s(t_f)\\
       x_{6,i}-q(t_f) 
    \end{bmatrix} = \boldsymbol{0}.
\end{equation}
With reference to $\tau_B$, the final state corresponds to orbit injection into a circular LLO, with prescribed radius and inclination. This leads to defining the final conditions, in the form
\begin{equation}\label{Psi_f_b}
    \boldsymbol{\Psi}_f^{(b)}=  \begin{bmatrix}
       x_{1,f}-p_d \\ 
       x_{2,f}^2+x_{3,f}^2-e_d^2\\ 
       x_{4,f}^2+x_{5,f}^2-\tan^2{\cfrac{i_d}{2}}
    \end{bmatrix} = \boldsymbol{0}.
\end{equation}
For case \ref{LLOToGat}, because backward propagation is employed, the state equations (\ref{dotx}) must be rewritten with $\tau_B$ as the independent variable, i.e.,
\begin{equation}\label{stateEqback}
     \dot{\boldsymbol{x}} =\boldsymbol{f} \left( \boldsymbol{x}, \boldsymbol{u}, \boldsymbol{p}, \tau_B \right) = -\boldsymbol{\tilde{f}} \left( \boldsymbol{x}, \boldsymbol{u}, t_f-\tau_B \right). 
\end{equation}
where the dot denotes the derivative with respect to $\tau_B$. Moreover, the thrust acceleration magnitude $a_T$ is given by
\begin{equation}\label{aTback}
     a_T=\frac{u_T^{(max)}}{m_R}=\frac{u_T^{(max)}}{1-\frac{u_T^{(max)}}{c}\left(\tau_{fin}-\tau_B\right)}
\end{equation}

For both cases (a) and (b), the objective function to minimize is the time of flight, i.e.,
\begin{equation}\label{objFun}
    J=k_J \tau_{fin}\,(=k_J(t_f-t_0))
\end{equation}
where either $\tau=\tau_F$ or $\tau=\tau_B$, while $k_J$ denotes an arbitrary positive constant.

\subsection{Method of solution}
The method of solution employs the indirect heuristic method \cite{pontani2014optimal,pontani2015minimum}, based on the joint use of the necessary conditions for an extremal, described in the Appendix, and a heuristic technique, i.e., differential evolution (DE) in this research.

As a preliminary step, based on previous research \cite{pozzi2024,beolchi2024}, the Hamiltonian function is rewritten as
\begin{equation}\label{Hsplit}
H=H_x+H_r a_{T,r}+H_{\theta} a_{T,\theta}+H_h a_{T,h}
\end{equation}
where $H_x$ denotes the portion of Hamiltonian that is independent of the control variables (associated with $a_{T,r}$, $a_{T,\theta}$, and $a_{T,h}$, cf. Eq. (\ref{thrustComp})), and 
\begin{align}
    H_{{r}} &= \sqrt{\frac{x_1}{\mu}}  \left( \lambda_{2}  \sin{x_6} - \lambda_{3}  \cos{x_6} \right)\\
    \begin{split}
        H_{{\theta}} &= \frac{\sqrt{\frac{x_1}{\mu}}}{1 + x_2  \cos{x_6} + x_3  \sin{x_6}}  \left\{ 2  x_1  \lambda_{1}  \right.\\
        &+ \left[x_3 + \cos{x_6}  \left(2 + x_2  \cos{x_6} + x_3  \sin{x_6} \right) \right]  \lambda_{2} \\
        &+ \left. \left[ x_3 + \sin{x_6}  \left(2 + x_2  \cos{x_6} + x_3  \sin{x_6} \right) \right]  \lambda_{3} \right\}
    \end{split}\\
    \begin{split}
        H_{{h}} &= \frac{\sqrt{\frac{x_1}{\mu}}}{2  \left( 1 + x_2  \cos{x_6} + x_3  \sin{x_6} \right)}  \left[  2  \left( x_4  \sin{x_6}  \right. \right.\\
        &- \left. x_5  \cos{x_6} \right)  \left( \lambda_{6} + x_2  \lambda_{3} - x_3  \lambda_{2} \right) \\
        &+ \left. \left( x_4^2 + x_5^2 + 1 \right)  \left( \cos{x_6}  \lambda_{4} + \sin{x_6}  \lambda_{5} \right) \right].
    \end{split}
\end{align}
 Using the relation in Eq. (\ref{Hsplit}), the Pontryagin minimum principle provides the control angles in terms of the state and costate variables \cite{pozzi2024,beolchi2024},
\begin{equation}
  \begin{aligned}\label{anglesalfa}
        \sin{u_2} &= -\frac{H_r}{\sqrt{H_r^2+H_{\theta}^2}} \, \, \, \, \, \, 
        \cos{u_2} = -\frac{H_\theta}{\sqrt{H_r^2+H_{\theta}^2}}
  \end{aligned}
\end{equation}
\begin{equation}\label{anglesbeta}
 \sin{u_3} = -\frac{H_h}{\sqrt{H_r^2+H_{\theta}^2+H_h^2}}\\
\end{equation}
Moreover, due to the scalability property of the adjoint variables \cite{pontanibook,beolchi2024}, their values at $\tau_i$, $\boldsymbol{\lambda}_i$, can be sought in the interval $\left[-1, \, \, \, 1\right]$. 

With regard to problem (a) (Gateway-to-LLO), the boundary conditions (\ref{Psi_0_a})-(\ref{Psi_f_a}), together with the necessary conditions reported in the Appendix, yield
\begin{equation}\label{lambdainia}
   \boldsymbol{\lambda}_i=-\boldsymbol{\upsilon}_i
\end{equation}
\begin{equation}\label{lambdafin1a}
   \lambda_{2,f} x_{3,f} -\lambda_{3,f} x_{2,f} =0
\end{equation}
\begin{equation}\label{lambdafin2a}
   \lambda_{4,f} x_{5,f} -\lambda_{5,f} x_{4,f} =0
\end{equation}
\begin{equation}\label{lambdafin3a}
   \lambda_{6,f}=0. 
\end{equation}
Moreover, because $H$, $\boldsymbol{\Psi}_i$, and $\boldsymbol{\Psi}_f$ 
are independent of $p_2\,(=t_f)$, Eq. (\ref{paramCond2}) yields $\xi=0$. The Two-Point Boundary-Value Problem (TPBVP) associated with the minimum-time problem involves the following ($8 \times 1$)-vector of unknown quantities:
\begin{equation}\label{XvecF}
    \boldsymbol{X} = \begin{bmatrix}
      p_1 & \tau_{fin} & \lambda_{1,i} & \lambda_{2,i} & \lambda_{3,i} & \lambda_{4,i} & \lambda_{5,i} & \lambda_{6,i} \\ 
    \end{bmatrix}^T. 
\end{equation}
\noindent The corresponding ($7 \times 1$) constraint vector includes the final boundary conditions on the state (\ref{Psi_f_a}) and the costate (\ref{lambdafin1a})-(\ref{lambdafin3a}), together with the parameter condition (\ref{paramCond1}), 
\begin{equation}\label{equalirycontrF}
    \boldsymbol{Y} = \begin{bmatrix}
       x_{1,f}-p_d \\ 
       x_{2,f}^2+x_{3,f}^2-e_d^2\\
       x_{4,f}^2+x_{5,f}^2-\tan^2{\cfrac{i_d}{2}} \\ \lambda_{2,f}{x_{3,f}} - \lambda_{3,f}{x_{2,f}} \\  \lambda_{4,f}{x_{5,f}} - \lambda_{5,f}{x_{4,f}} \\ \lambda_{6,f} \\
       \int_{\tau_i}^{\tau_{fin}}\left(\frac{\partial H}{\partial p_1}\right) d\tau
       -\left(\frac{\partial\boldsymbol{\Psi_i}}{\partial p_1}\right)^T\boldsymbol{\lambda}_i
    \end{bmatrix} = \boldsymbol{0},
\end{equation}
\noindent where $\left(\partial\boldsymbol{\Psi_i}/\partial p_2=\boldsymbol{0}\right)$ and Eq. (\ref{lambdainia}) were used in the last relation. The remaining constraint is represented by the inequality condition on the final Hamiltonian (\ref{Hf}), with $\xi=0$.

With regard to problem (b) (LLO-to-Gateway), the boundary conditions (\ref{Psi_0_b})-(\ref{Psi_f_b}), together with the necessary conditions reported in the Appendix, yield exactly the same relations (\ref{lambdainia})-(\ref{lambdafin3a}). Moreover, because $H$, $\boldsymbol{\Psi}_i$, and $\boldsymbol{\Psi}_f$ 
are independent of $p_1\,(=t_0)$, Eq. (\ref{paramCond1}) yields $\xi=0$. 
The Two-Point Boundary-Value Problem (TPBVP) associated with the minimum-time problem at hand involves the following ($8 \times 1$)-vector of unknown quantities:
\begin{equation}\label{XvecB}
    \boldsymbol{X} = \begin{bmatrix}
      p_2 & \tau_{fin} & \lambda_{1,i} & \lambda_{2,i} & \lambda_{3,i} & \lambda_{4,i} & \lambda_{5,i} & \lambda_{6,i} \\ 
    \end{bmatrix}^T. 
\end{equation}
\noindent The corresponding ($7 \times 1$) constraint vector includes the final boundary conditions on the state (\ref{Psi_f_a}) and the costate (\ref{lambdafin1a})-(\ref{lambdafin3a}), together with the parameter condition (\ref{paramCond2}), 
\begin{equation}\label{equalirycontrB}
    \boldsymbol{Y} = \begin{bmatrix}
       x_{1,f}-p_d \\ 
       x_{2,f}^2+x_{3,f}^2-e_d^2\\
       x_{4,f}^2+x_{5,f}^2-\tan^2{\cfrac{i_d}{2}} \\ \lambda_{2,f}{x_{3,f}} - \lambda_{3,f}{x_{2,f}} \\  \lambda_{4,f}{x_{5,f}} - \lambda_{5,f}{x_{4,f}} \\ \lambda_{6,f} \\
       \int_{\tau_i}^{\tau_{fin}}\left(\frac{\partial H}{\partial p_2}\right) d\tau-
       \left(\frac{\partial\boldsymbol{\Psi_i}}{\partial p_2}\right)^T\boldsymbol{\lambda}_i
    \end{bmatrix} = \boldsymbol{0},
\end{equation}
\noindent where $\left(\partial\boldsymbol{\Psi_i}/\partial p_1=\boldsymbol{0}\right)$ and Eq. (\ref{lambdainia}) were used in the last relation. The remaining constraint is represented by the inequality condition on the final Hamiltonian (\ref{Hf}), with $\xi$ given by Eq. (\ref{paramCond1}). It is worth remarking that $p_1$ can be evaluated using the relation $p_1=p_2-\tau_{fin}$. 

In Appendix B, the last two conditions that appear in the two distinct vectors $\boldsymbol{Y}$ in Eqs. (\ref{equalirycontrF}) and (\ref{equalirycontrB}) can be proven to be equivalent to two inequality constraints, and turn out to be always satisfied at the optimal solution, as a consequence of the Pontryagin minimum principle. As a result, the two scalar equality constraints associated with the last two elements of $\boldsymbol{Y}$ can be neglected in the numerical solution process.

As a preliminary step for the numerical solution process, canonical units are introduced. The Distance Unit ($\rm DU$) equals the mean equatorial radius of the main attracting body (the Moon), whereas the Time Unit ($\rm TU$) is such that the lunar gravitational parameter is  $\mu = 1 \, \, {\rm DU^3}/{\rm TU^2}$. 
 The numerical solution process utilizes a population of individuals, each corresponding to a possible solution. The fitness function is defined in the next subsection. Based on the necessary conditions, the solution technique consists of the following steps: 

 \begin{enumerate}[itemsep=-0.1cm]
     \item Identify the known initial values of the state and costate variables and the minimal set of unknown values.
     \item For each individual $i$, with $i \leq N_p$ (with $N_p$ denoting the total number of individuals), iterate the following sub-steps:
\begin{enumerate}
    \item select the unknown values for the quantities contained in $\boldsymbol{Y}$ (cf. Eqs. (\ref{equalirycontrF}) and (\ref{equalirycontrB})); 
    \item after selecting either $t_0$ or $t_f$, calculate the initial values of the state components;
    \item integrate numerically the state and costate equations while using Eqs. (\ref{anglesalfa})-(\ref{anglesbeta}) until $\tau = \tau_{fin}$;
    \item evaluate the Hamiltonian at the final time and the integral term that appears in the last component of $\boldsymbol{Y}$ (cf. Eqs. (\ref{equalirycontrF}) and (\ref{equalirycontrB}));
    \item if inequality (\ref{Hf}) is met, then go to sub-step (f), otherwise set the fitness function $\tilde{J}$ to a large value and go on to the next individual;
    \item evaluate the violations of the boundary conditions (cf. Eqs (\ref{equalirycontrF}) or (\ref{equalirycontrB})), then compute the fitness function $\tilde{J}$, defined in the next subsection;
    \item go to the next individual.
\end{enumerate}
     \item Once the fitness functions are obtained for all the individuals that compose the population, use the DE algorithm to update the position of individuals in the search space.
     \item Repeat all the steps until the fitness function of the best individual reaches a value lower than a threshold  or until the maximum number of iterations is reached.
     
 \end{enumerate}

\indent DE stops when a prescribed number of iterations is reached or if the global best has not changed for several stalled iterations. In this work, the same settings adopted in \cite{beolchi2024} are used. 

\subsection{Numerical results}
\indent The point-mass spacecraft is characterized by the following propulsion parameters:
\begin{equation}\label{proppar}
        u_T^{(max)}  = 4.903 \cdot 10^{-7} \, \, \rm{\cfrac{km}{s^2}} \quad \textrm{and}\quad 
        c = 30 \, \, \rm{\cfrac{km}{s}}.
\end{equation}
If the space vehicle has mass of 15000 kg (e.g., ATV \cite{amadieu1999automated}), then the maximum thrust magnitude equals $T^{(max)}=7.3545 \cdot 10^{-3}$ N.

The objective of the optimization algorithm is an equality-constraint-inspired fitness function $\Tilde{J}$ defined as
\begin{equation}\label{Jtilde}
    \tilde{J}=\sqrt{\sum_{i=1}^{7}w_i Y_i^2}
\end{equation}
\noindent where $Y_i$ (for $i = 1 \dots 7$) are the elements of the constraint vector in either Eq. (\ref{equalirycontrF}) or Eq. (\ref{equalirycontrB}), and the $w_i$ are weights specifically chosen to give comparable relevance to the boundary conditions to satisfy. After a tuning phase, the  weights are
\begin{equation}
    w_1=w_3=w_4=w_5=w_6=w_7 = 1 \, \, \,w_2 = 100.
\end{equation}
The choice of setting $w_2$ to 100 is justified by the fact that convergence for the eccentricity turned out to be more challenging. The initial values of the adjoint variables are always sought in the interval $[-1,1]$, due to their scalability, while the time of flight is sought in the interval [25,45] days, based on the analytical estimate described in \cite{pozzi2024}.

For case (a) (Gateway-to-LLO transfer), the spacecraft of interest travels the NRHO flown by Gateway, and the initial time of the transfer (i.e., $t_0$) is selected along a single period of the departing orbit, in the interval [23 May 2025 at 22:35, 30 May 2025 at 05:38] UTC. The final orbit is a LLO with the following desired orbit elements: 
\begin{equation}\label{edid}
p_d = R_{M} + 100 \, \, {\rm{km}} \, \, \, \, \,  e_d = 0 \, \, \, \,  \,  i_d = 90^{\circ}
\end{equation}
where $R_{M} = 1737.4 \, \, \rm{km}$ is the mean equatorial radius of the Moon.

\begin{figure}[H]
\centering
   \includegraphics[width=0.92\columnwidth]{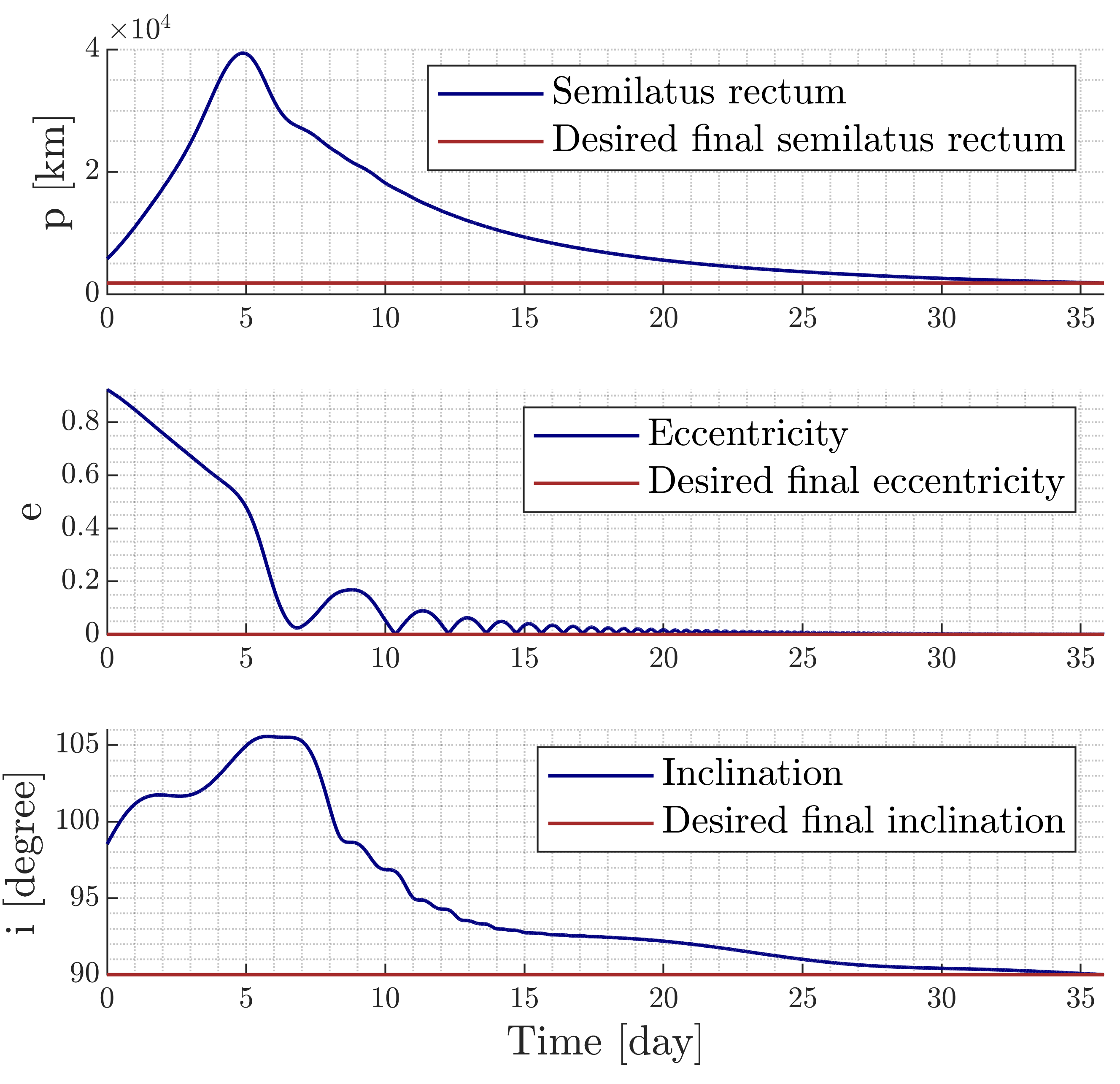}
   \caption{Gateway-to-LLO: time histories of the relevant orbit elements}
    \label{pei_optF}
\end{figure}

\begin{figure}[H]
     \includegraphics[width=\columnwidth]{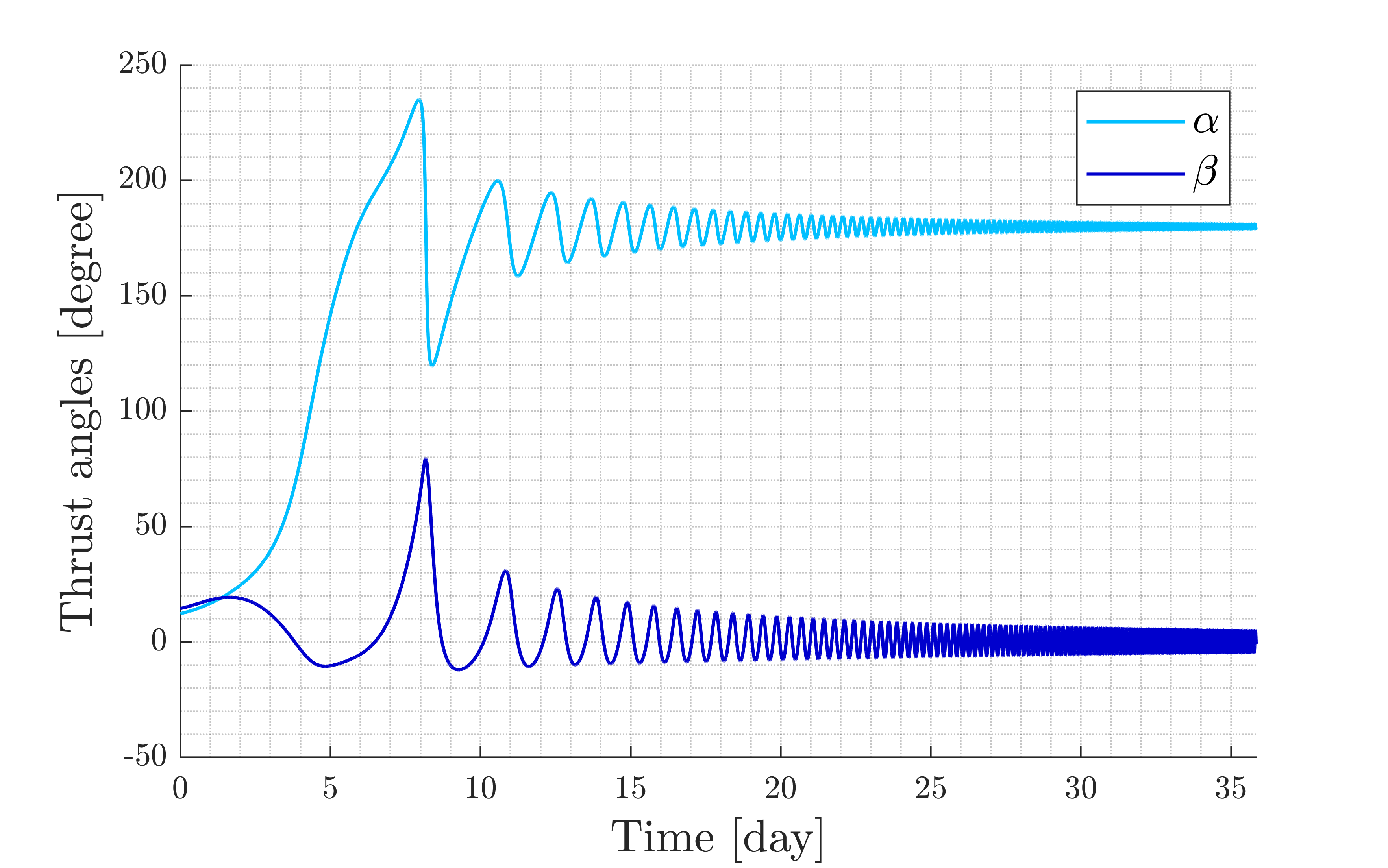}
    \caption{Gateway-to-LLO: time histories of thrust angles}
    \label{alfabeta_optF}
\end{figure}

\begin{figure}[H]
\centering
    \includegraphics[width=0.95\columnwidth]{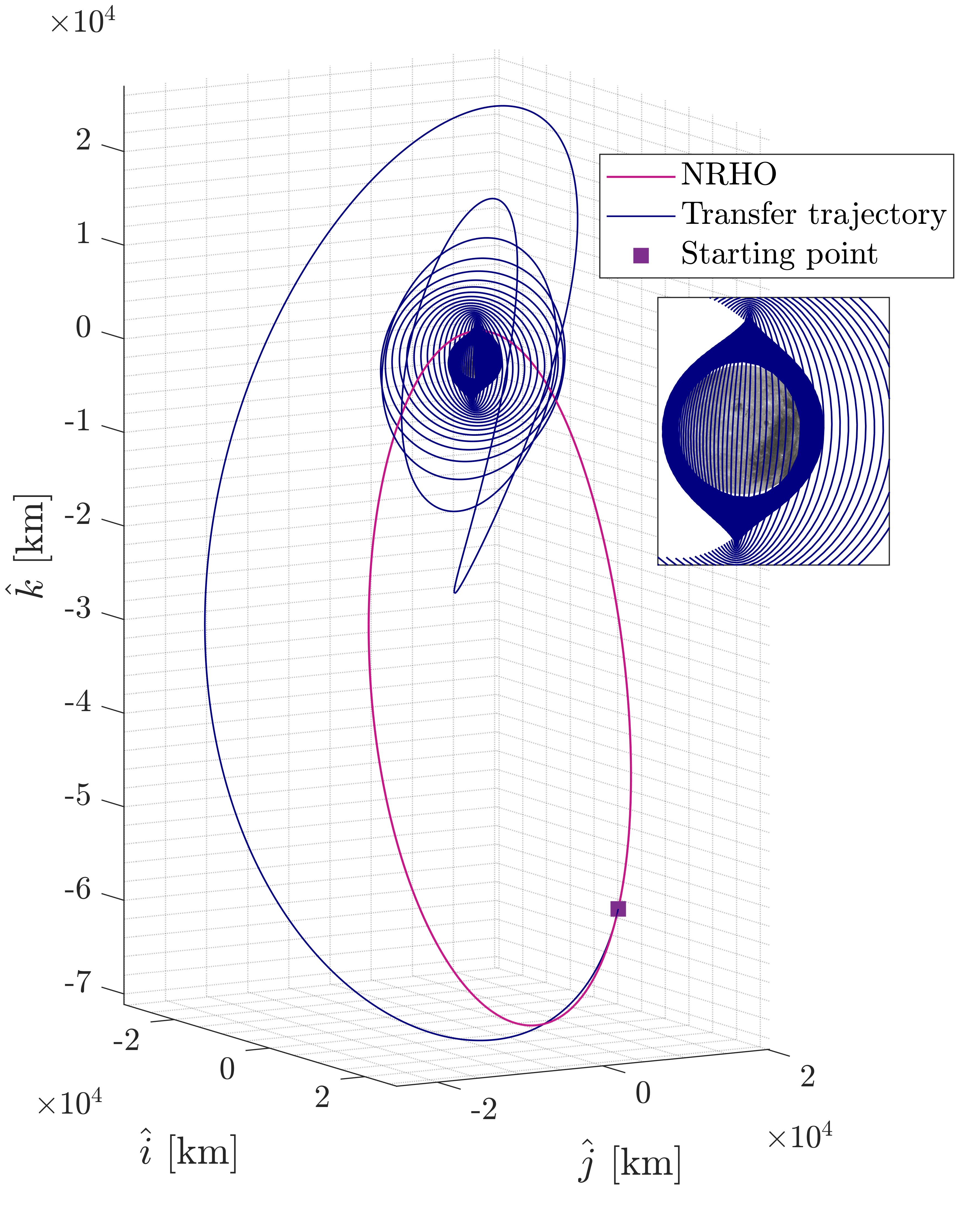}
    \caption{Gateway-to-LLO: minimum-time transfer in the synodic frame}
    \label{trajoptF}
\end{figure}

\indent After running DE, the \textit{fminsearch} MATLAB routine is employed for refinement, leading to $\Tilde{J} = 1.911 \cdot 10^{-10}$. 
The optimal starting time is 25 May 2025 at 16:51:30 UTC and the time of flight equals $ \tau_{fin} = 35 \, \, {\rm d} \, \, 20 \, \, {\rm hrs} \, \, 2 \, \, {\rm min} \, \, 45 \, \, {\rm s} $, whereas the final mass ratio is $0.9494$. Table \ref{tabstate} provides the spacecraft orbit elements at the initial and final time. The orbit elements at $t_0$ are the osculating elements of Gateway at the beginning of the transfer. Figure \ref{pei_optF} depicts the time histories of $p$, $e$, and $i$, whereas Figs. \ref{alfabeta_optF} and \ref{trajoptF} portray the optimal thrust angles and transfer trajectory. The time history of the eccentricity shows a rapid decrease in the first 6 days, then it reduces to 0 in the remaining days. From inspection of the time evolution of the inclination, it is apparent that the orbital plane changes continuously, although modest variations occur after 15 days. Moreover, the thrust direction is nearly opposed to the spacecraft velocity in the last 20 days, to reduce the orbit eccentricity and semilatus rectum. Figure \ref{trajoptF} represents the optimal transfer trajectory in the Moon-centered synodic frame. The starting point is close to aposelenium of the NRHO and the spacecraft completes many orbits about the Moon before arriving at the final desired LLO.

\begin{table}[H] 
\centering
\begin{tabular}{ccc}
\toprule
\text{COE} & ${t_0}$ & ${t_f}$\\
\midrule
$a$ [km]&$3.916\cdot10^{4}$&$1.837\cdot10^{3}$\\ 
$e$&$0.923$ & $1.289\cdot10^{-6}$\\ 
$i\,[\degree]$  &$98.53 $ & $90.00 $\\
$\Omega\,[\degree]$  &$ -60.75$ &$-47.86$ \\
$\omega\,[\degree]$  &$84.05 $ & $28.08$ \\
$\theta_*\,[\degree]$  &$168.22$ & $-109.91$  \\
\bottomrule
\end{tabular}
\vspace*{2mm}
\caption{Gateway-to-LLO: initial and final orbit elements}
\label{tabstate}
\end{table}

For case (b) (LLO-to-Gateway transfer), initially the spacecraft travels a circular, polar LLO and the final time of the transfer (i.e., $t_f$) is selected in the same time interval of case (a), referred to the period of the NRHO traveled by Gateway. More precisely, the orbit at $t_0$ is a LLO with the following orbit elements: 
\begin{equation}\label{edid}
p_d = R_{M} + 100 \, \, {\rm{km}} \, \, \, \, \,  e_d = 0 \, \, \, \,  \,  i_d = 90^{\circ}
\end{equation}

\indent After running DE, the \textit{fminsearch} MATLAB routine is employed for refinement, leading to $\Tilde{J} =7.201 \cdot 10^{-8}$. 
The optimal starting time is 22 April 2025 at 16:39:52 UTC and the time of flight equals $ \tau_{fin}= 36 \, \, {\rm d} \, \, 4\, \, {\rm hrs} \, \, 24 \, \, {\rm min} \, \, 52 \, \, {\rm s} $, whereas the final mass ratio is $0.9489$. Table \ref{tabstate2} provides the spacecraft orbit elements at the initial and final time. The orbit elements at $t_f$ correspond to the osculating elements of Gateway at $t_f$.  

\bgroup
\def\arraystretch{1.1}
\begin{table}[h] 
\centering
\begin{tabular}{ccc}
\toprule
\text{COE} & ${t_0}$ & ${t_f}$\\
\midrule
$a$ [km]&$1.837\cdot10^{3}$&$3.931\cdot10^{4}$\\ 
$e$&$2.112\cdot10^{-5}$ & $0.923$\\ 
$i\,[\degree]$  &$90.00$ & $100.99$\\
$\Omega\,[\degree]$  &$-43.18$ &$ 11.10$\\
$\omega\,[\degree]$  &$-87.63$ & $92.18$\\
$\theta_*\,[\degree]$  &$27.30$ & $-167.67$\\
\bottomrule
\end{tabular}
\vspace*{2mm}
\caption{LLO-to-Gateway: initial and final orbit elements}
\label{tabstate2}
\end{table}
\egroup

\begin{figure}[H]
   \includegraphics[width=1\columnwidth]{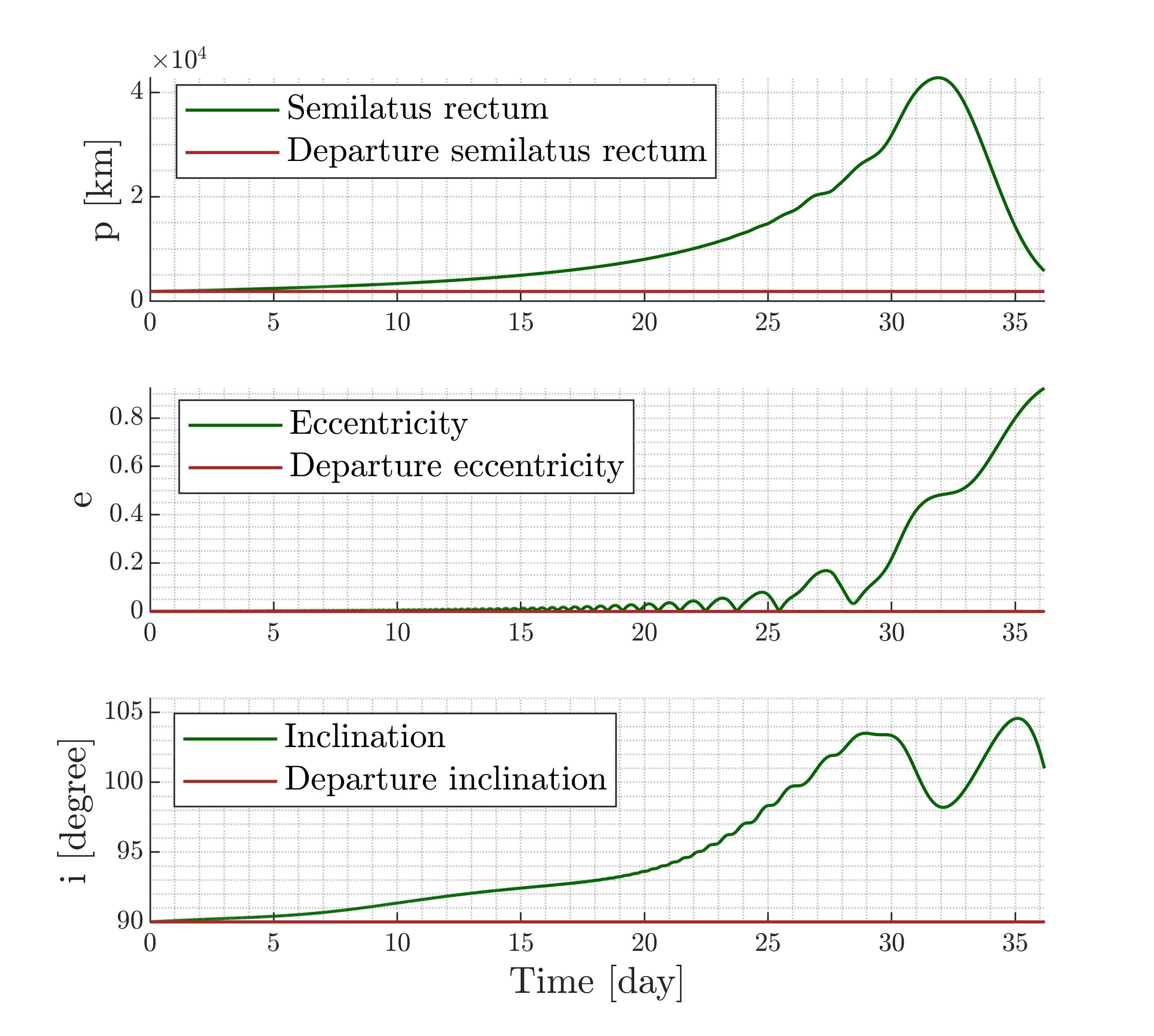}
   \caption{Gateway-to-LLO: time histories of the relevant orbit elements}
    \label{pei_optB}
\end{figure}

\begin{figure}[H]
     \includegraphics[width=\columnwidth]{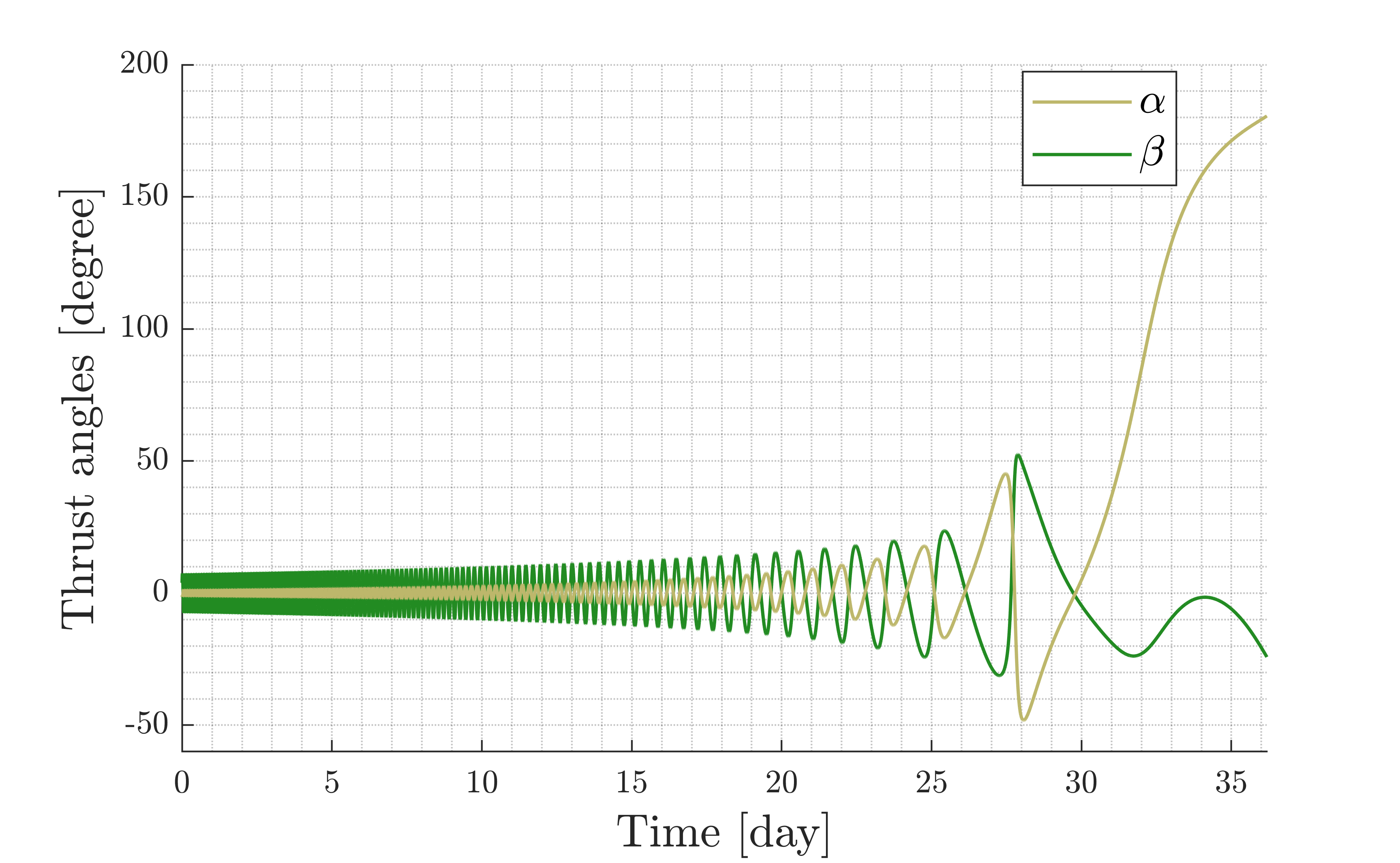}
    \caption{LLO-to-Gateway: time histories of thrust angles}
    \label{alfabeta_optB}
\end{figure}

\begin{figure}[H]
    \includegraphics[width=\columnwidth]{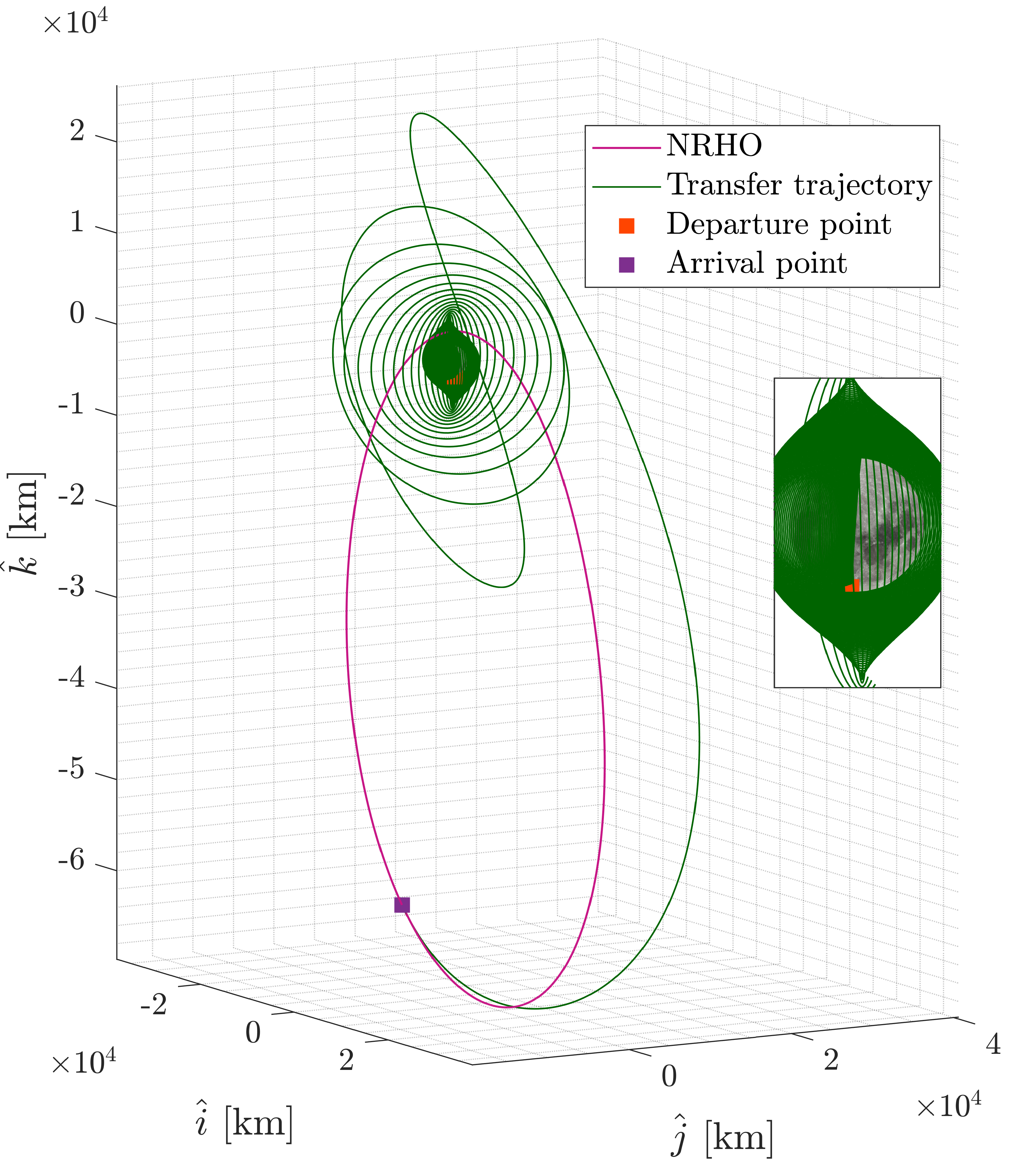}
    \caption{LLO-to-Gateway: minimum-time transfer in the synodic frame}
    \label{trajoptB}
\end{figure}

Figure \ref{pei_optB} depicts the time histories of $p$, $e$, and $i$, whereas Figs. \ref{alfabeta_optB} and \ref{trajoptB} portray the optimal thrust angles and transfer trajectory. The time history of the eccentricity shows a very slow increase in the first 29 days, then it rises to the desired final value in the last 6 days. From inspection of the time evolution of the inclination it is apparent that the orbital plane changes continuously, although modest variations occur in the first 20 days. Moreover, the thrust direction is aligned with the spacecraft velocity in the first 20 days, to increase eccentricity and semilatus rectum. Figure \ref{trajoptB} represents the optimal transfer trajectory in the Moon-centered synodic frame. The arrival point is close to aposelenium of NRHO.



\section{LOW-THRUST TRANSFERS CONNECTING GATEWAY AND LEO} 

This section is focused on identifying the minimum-time transfer paths that connect Gateway and a LEO. In particular, two distinct low-thrust orbit transfers are considered: 
\begin{enumerate}[label=(\small{\alph*}),itemsep=0cm]
\vspace{-0.2cm}\item\label{GatToLEO} transfer from Gateway to LEO, and
\vspace{-0.2cm}\item\label{LEOToGat} transfer from LEO to Gateway.
\end{enumerate}
Unlike the previous case, the problem at hand requires using both terrestrial and lunar MEE, because the terminal orbits have two different dominating attracting bodies. This implies the need of identifying suitable transformations for the time-varying quantities involved in the solution process, primarily the state and costate variables.

Previous research also proved that for minimum-time transfers the thrust magnitude must be set to its maximum value, thus Eq. (\ref{u_1}) holds. This circumstance allows integrating $x_7$ in closed form. As a result, this variable can be eliminated from the state vector.

\subsection{Formulation of the problem}
\indent In case \ref{GatToLEO}, forward propagation is used, by introducing the auxiliary forward time variable $\tau_F:=t-t_0$. Both the initial time $t_0$ and the final time $t_f$ are unspecified, and $\tau_F$, which plays the role of independent variable for the dynamical system at hand, is constrained to $[0,t_f-t_0]$. Time $t_0$ identifies the initial state (cf. Eq. (\ref{stateMEE})) at departure, which is collected in $\boldsymbol{\Psi}_i^{(a)}$, which fulfills 
\begin{equation}\label{Psi_0_a_2}
    \boldsymbol{\Psi}_i^{(a)}=  \begin{bmatrix}
       x_{M,1,i}-p_M(t_0) \\
       x_{M,2,i}-l_M(t_0)\\
       x_{M,3,i}-m_M(t_0)\\
       x_{M,4,i}-n_M(t_0)\\
       x_{M,5,i}-s_M(t_0)\\
       x_{M,6,i}-q_M(t_0)        
    \end{bmatrix} = \boldsymbol{0}.
\end{equation}
where subscript $M$ refers to the Moon (i.e., in Eq. (\ref{Psi_0_a_2}) lunar MEE appear). The final state corresponds to orbit injection into a circular LEO, with prescribed radius and inclination. This leads to defining the final conditions, in the form
\begin{equation}\label{Psi_f_a_2}
    \boldsymbol{\Psi}_f^{(a)}=  \begin{bmatrix}
       x_{E,1,f}-p_{E,d} \\ 
       x_{E,2,f}^2+x_{E,3,f}^2-e_{E,d}^2\\ 
       x_{E,2,f}^2+x_{E,5,f}^2-\tan^2{\cfrac{i_{E,d}}{2}}
    \end{bmatrix} = \boldsymbol{0}
\end{equation}
where subscript $E$ refers to Earth (i.e., in Eq. (\ref{Psi_f_a_2}) terrestrial MEE appear), whereas $p_{E,d}=R_{LEO}$ denotes the desired semilatus rectum, $e_{E,d}$ $(=0)$ the final eccentricity, and $i_{E,d}$ the desired inclination of LEO. For case \ref{GatToLEO}, because forward propagation is employed, the state equations (\ref{dotx}) remain unaltered when rewritten with $\tau_F$ as the independent variable, i.e,
\begin{equation}\label{stateEqforw2}
     \dot{\boldsymbol{x}} =\boldsymbol{f} \left( \boldsymbol{x}, \boldsymbol{u}, \boldsymbol{p}, \tau_F \right) = \boldsymbol{\tilde{f}} \left( \boldsymbol{x}, \boldsymbol{u}, \tau_F+t_0 \right). 
\end{equation}
where the dot denotes the derivative with respect to $\tau_F$, whereas $\boldsymbol{p}$ represents a vector of unknown time-independent parameters, i.e., the two epochs $t_0$ and $t_f$ in this research. Moreover, the thrust acceleration magnitude $a_T$ is given by Eq. (\ref{aTforw}). 

\indent In case \ref{LEOToGat}, the spacecraft travels from a circular LEO with specified radius $R_{LEO}$ and inclination $i_{E,d}$ toward Gateway. Because the final dynamical conditions, i.e., rendezvous with Gateway, are more stringent than the initial conditions, it is convenient to investigate this transfer using backward time propagation. To this end, the auxiliary backward time variable $\tau_B:=t_f-t$ is introduced. Both the initial (actual) time $t_0$ and the final actual time $t_f$ are unspecified, therefore $\tau_B$, which plays the role of independent variable for the dynamical system at hand, is monotonically increasing and constrained to $[0,t_f-t_0]$. With respect to $\tau_B$, the actual time $t_f$ identifies the initial state, which fulfills 
\begin{equation}\label{Psi_0_b_2}
    \boldsymbol{\Psi}_i^{(b)}=  \begin{bmatrix}
       x_{M,1,i}-p_M(t_f) \\
       x_{M,2,i}-l_M(t_f)\\
       x_{M,3,i}-m_M(t_f)\\
       x_{M,4,i}-n_M(t_f)\\
       x_{M,5,i}-s_M(t_f)\\
       x_{M,6,i}-q_M(t_f)        
    \end{bmatrix} = \boldsymbol{0}.
\end{equation}
With reference to $\tau_B$, the final state corresponds to orbit injection into a circular LEO, with prescribed radius and inclination. This leads to defining the final conditions, in the form
\begin{equation}\label{Psi_f_b_2}
    \boldsymbol{\Psi}_f^{(b)}=  \begin{bmatrix}
       x_{E,1,f}-p_{E,d} \\ 
       x_{E,2,f}^2+x_{E,3,f}^2-e_{E,d}^2\\ 
       x_{E,4,f}^2+x_{E,5,f}^2-\tan^2{\cfrac{i_{E,d}}{2}}
    \end{bmatrix} = \boldsymbol{0}.
\end{equation}
For case \ref{LEOToGat}, because backward propagation is employed, the state equations (\ref{dotx}) must be rewritten with $\tau_B$ as the independent variable, i.e.,
\begin{equation}\label{stateEqback2}
     \dot{\boldsymbol{x}} =\boldsymbol{f} \left( \boldsymbol{x}, \boldsymbol{u}, \boldsymbol{p}, \tau_B \right) = -\boldsymbol{\tilde{f}} \left( \boldsymbol{x}, \boldsymbol{u}, t_f-\tau_B \right). 
\end{equation}
where the dot denotes the derivative with respect to $\tau_B$. Moreover, the thrust acceleration magnitude $a_T$ is given by Eq. (\ref{aTback}).

It is worth remarking that the preceding state equations (\ref{stateEqforw2}) and (\ref{stateEqback2}) are formally identical for terrestrial and lunar MEE. These two different subsets are being distinguished through superscripts in next sections.

For both cases (a) and (b), the objective function to minimize is the time of flight, i.e.,
\begin{equation}\label{objFun}
    J=k_J \tau_{fin}\,(=k_J(t_f-t_0))
\end{equation}
where either $\tau=\tau_F$ or $\tau=\tau_B$, while $k_J$ denotes an arbitrary positive constant.

\subsection{Trajectory arcs and coordinate transformations}\label{TACT} 

\indent In both cases (a) and (b) terrestrial MEE are more convenient in the proximity of Earth, whereas lunar MEE are used while approaching to (or departing from) Gateway.
Thus, the transfer trajectory can be virtually partitioned in two arcs: (i) a geocentric arc and (ii) a selenocentric arc. The timing for transitioning from (i) to (ii) (or vice-versa) is arbitrary \cite{beolchi2024}, but proper selection of the transition time facilitates the convergence of the numerical solution process. In this study, the change in representation occurs as soon as the spacecraft has a distance from Earth equal to 320000 km, denoted with $\rho_E$, i.e., when the following equality is satisfied: 
\begin{equation}\label{SOImoon}
    r_{E} - \rho_E = 0
\end{equation} 
\noindent where $r_{E}$ represents the instantaneous distance of the space vehicle from the Earth. The orbit of Gateway is completely out of the sphere centered at Earth and with radius $\rho_E$, and this allows describing the Gateway orbital motion with the use of selenocentric MEE.

For both cases (a) and (b), 5 distinct state representations are needed, i.e., (1) lunar MEE $\boldsymbol{x}_M$, (2) lunar CC $\boldsymbol{y}_M$, (3) lunar CC projected onto the ECI-frame $\boldsymbol{y}_M^{(ECI)}$, (4) terrestrial CC $\boldsymbol{y}_E$, and (5) terrestrial MEE $\boldsymbol{x}_E$. Moreover, in both cases the initial time in the $\tau$-domain corresponds to departure from Gateway, either with forward propagation (case (a), Gateway-to-LEO transfer) or with backward propagation (case (b), LEO-to-Gateway transfer).

\indent With reference to both cases (a) and (b), as $\tau$ increases $\boldsymbol{x}_E$ is obtained from $\boldsymbol{x}_M$ through 4 transformations:
\begin{enumerate}[label={(\Alph*)}]
    \item\label{transf1a} $\boldsymbol{x}^{(1)}:=\boldsymbol{x}_M$ to $\boldsymbol{x}^{(2)}:=\boldsymbol{y}_M$
    \item\label{transf2a} $\boldsymbol{x}^{(2)}:=\boldsymbol{y}_M$ to $\boldsymbol{x}^{(3)}:=\boldsymbol{y}_M^{(ECI)}$    
    \item\label{transf3a} $\boldsymbol{x}^{(3)}:=\boldsymbol{y}_M^{(ECI)}$ to $\boldsymbol{x}^{(4)}:=\boldsymbol{y}_E$    
    \item\label{transf4a} $\boldsymbol{x}^{(4)}:=\boldsymbol{y}_E$ to $\boldsymbol{x}^{(5)}:=\boldsymbol{x}_E$
\end{enumerate}
While steps \ref{transf1a} and \ref{transf4a} represent nonlinear isomorphic (bijective) mappings between CC (of position and velocity) and MEE, step \ref{transf2a} is a rotation. Finally, \ref{transf3a} represents a translation of the spacecraft position and velocity vectors, which refer to two distinct bodies, in the sense that $\boldsymbol{y}_M^{(ECI)}$ and $\boldsymbol{y}_E$ collect the components relative to Moon and Earth, respectively (both projected onto the ECI-frame). It is worth noting that state transformations \ref{transf1a}, \ref{transf2a}, and \ref{transf4a} do not explicitly depend on time.

The $j$-th step of the previous sequence, occurring at time $\tau_j$, can be represented through the general expression for implicit state transformations
\begin{equation} \label{chij}
    \boldsymbol{\chi}_j \left( \boldsymbol{x}_{ini}^{(j + 1)}, \boldsymbol{x}_{fin}^{(j)}, \tau_j \right) = \boldsymbol 0
\end{equation}
\noindent where $\boldsymbol{x}_{ini}^{(j + 1)}$ and $\boldsymbol{x}_{fin}^{(j)}$ denote the state after and before transformation $j$, respectively, and $\boldsymbol{\chi}_j$ is a nonlinear vector function, whose dimension equals the number of state components. Hence forward, subscripts $ini$ and $fin$ denote respectively the initial and final value of a variable in the arc with index reported in the superscript. Each transformation can also be expressed in explicit form, which guarantees a unique solution for $\boldsymbol{x}^{(j+ 1)}$, 
\begin{equation}\label{StateExplicitTransfTime} 
    \boldsymbol{x}_{ini}^{(j+ 1)} = \boldsymbol{\Lambda}_j \left( \boldsymbol{x}_{fin}^{(j)}, \tau_j \right)
\end{equation}

\indent For the problem at hand, even though the trajectory design strategy requires passing through three additional representations for the dynamical state of the spacecraft, only terrestrial and lunar MEE are used for the numerical propagation of the path of the space vehicle. This requires completing all four coordinate transformations simultaneously, i.e., $\tau_1 = \tau_2 = \tau_3 = \tau_4$.


\subsection{Multi-arc optimal control problem}\label{MultiarcOptContr} 
\indent For the orbit transfer problem at hand, the set of coordinates that represent the dynamical state of the spacecraft changes as soon as the space vehicle meets equality (\ref{SOImoon}). Because the entire trajectory is composed of five arcs (three of them have zero length), this becomes a multi-arc trajectory optimization problem, with number of arcs $N = 5$. In each arc $j$, the state equations are
\begin{equation}\label{dotxj}
    \dot{\boldsymbol{x}}^{(j)} = \boldsymbol{f}^{(j)} \left( \boldsymbol{x}^{(j)}, \boldsymbol{u}^{(j)},\boldsymbol{p},\tau \right) .
\end{equation}
The auxiliary independent variable $\tau$  plays the role of either forward time (for transfers that originate from Gateway) or backward time (for transfers that rendezvous with Gateway). In both cases, $\tau$ is monotonically increasing and constrained to $[0,t_f-t_0]$. The two epochs $t_0$ and $t_f$ are collected as components of the parameter vector $\boldsymbol{p}:=\left[p_1\,\,p_2\right]^T=\left[t_0\,\,t_f\right]^T$, and the following equality constraint is introduced:
\begin{equation}\label{taufinAndp_1}
\tau_{fin}+p_1-p_2=0        
    \end{equation}
\indent Denoting with $\tau_j$ the time at the $j$-th interface (between arc $j$ and arc $(j + 1)$), three additional functions are introduced:
\begin{itemize}[itemsep=-0.1cm]
    \item scalar transition function $\zeta_j$, identifying the occurrence of the transition between two consecutive arcs. In general, it is a function of $\boldsymbol{x}_{ini}^{(j + 1)}$, $\boldsymbol{x}_{fin}^{(j)}$, and $\tau_j$,
    \begin{equation}\label{zetaj} 
        {\zeta}_j \left( \boldsymbol{x}_{ini}^{(j + 1)}, \boldsymbol{x}_{fin}^{(j)}, \tau_j \right) = 0 \; 
    \end{equation}
    \item vector matching function $\boldsymbol{\chi}_j$, stating the (generally implicit) matching relation for the state across two adjacent arcs, i.e., Eq. (\ref{chij}),
    \item equality constraint (\ref{taufinAndp_1}).
\end{itemize}
\noindent The number of transitions is $(N - 1)$. 

\indent Optimal control theory is applied to the previously defined continuous-time dynamical system, for the purpose of obtaining the first-order necessary conditions for optimality, leading to translating the optimal control problem into a two-point boundary-value problem (TPBVP).

\indent As a first step, a function of boundary conditions $\Phi$ and $N$ Hamiltonian functions $H^{(j)}$, are introduced, 
\begin{equation}\label{Phimultiarc}
\begin{split}
    \Phi& = k_J \, \tau_{fin} + \boldsymbol{\upsilon}_i^T\boldsymbol{\Psi_i}+\boldsymbol{\upsilon}_f^T\boldsymbol{\Psi_f} \\
    &+ \sum\limits_{j = 1}^{N - 1} \left( \varepsilon_j \zeta_j + \boldsymbol{\nu}_j^{T} \boldsymbol{\chi}_j \right)+
    \xi\left(\tau_{fin}+p_1-p_2\right)
    \end{split}
\end{equation}
\begin{equation}\label{Hx}
    H^{(j)} = \boldsymbol{\lambda}^{{(j)}^T} \boldsymbol{f}^{(j)}
\end{equation}
\noindent where $N$ is the total number of arcs, $\boldsymbol{\upsilon}_i$, $\boldsymbol{\upsilon}_f$, $\varepsilon_j$, $\boldsymbol{\nu}_j$, and $\xi$ are time-independent adjoint variables conjugate respectively to the multipoint conditions (\ref{Psi_0_a_2}) and (\ref{Psi_f_a_2}) (or (\ref{Psi_0_b_2}) and (\ref{Psi_f_b_2})), (\ref{zetaj}), (\ref{chij}), and (\ref{taufinAndp_1}) , whereas $\boldsymbol{\lambda}^{(j)}$ is the time-varying costate vector associated with the differential constraint arising from the state equations (\ref{dotxj}). The extended objective functional for multi-arc problems is defined as
\begin{equation}\label{Jbarmultiarc}
    \begin{split}
        \overline{J} &= \Phi \left( \boldsymbol{x}_{i}, \boldsymbol{x}_{f}, \boldsymbol{p}, \boldsymbol{\sigma}, \tau_{fin}, \boldsymbol{\upsilon}_i, \boldsymbol{\upsilon}_f  \boldsymbol{\nu}_{1}, \dots,\boldsymbol{\nu}_{N - 1},\xi, \right. \\
        &  \boldsymbol{x}_{ini}^{(2)}, \dots, \boldsymbol{x}_{ini}^{(N)}, \boldsymbol{x}_{fin}^{(1)}, \dots, \boldsymbol{x}_{fin}^{(N - 1)}, \varepsilon_1, \dots,\left. \varepsilon_{N - 1} \right) \\ 
        &  + \sum\limits_{j = 1}^{N} \int_{\tau_{j - 1}}^{\tau_{j}} \left[ H^{(j)} \left( \boldsymbol{x}^{(j)}, \boldsymbol{u}^{(j)}, \boldsymbol{p}, \boldsymbol{\lambda}^{(j)}, \tau \right)  \right. \\
        & \left.- \boldsymbol{\lambda}^{{(j)}^T} \dot{\boldsymbol{x}}^{(j)} \right] d\tau
    \end{split}
\end{equation}
\noindent where $\tau_0 = \tau_i=0$, $\tau_N = \tau_{fin}$, while $\boldsymbol{x}_i:=\boldsymbol{x}_{ini}^{(1)}$ and $\boldsymbol{x}_f:=\boldsymbol{x}_{fin}^{(5)}$.

\indent The first differential of the augmented objective functional ${\rm d} \overline{J}$ can be obtained using the chain rule after lengthy {developments \cite{pontani2021optimal}}, omitted for the sake of brevity. Then, the necessary conditions for optimality are derived by enforcing ${\rm d} \overline{J} = 0$. They include the adjoint equations in each subarc and the Pontryagin minimum principle \cite{pontani2021optimal}, which leads to expressing the control angles in terms of costate variables (cf. Eq. (\ref{anglesalfa})). Moreover, Eqs. (\ref{lambda0}) and (\ref{lambdaf}) hold, where  $\boldsymbol{\lambda}_i$ and  $\boldsymbol{\lambda}_f$ are respectively identified with  $\boldsymbol{\lambda}_i^{(1)}$ and $\boldsymbol{\lambda}_f^{(5)}$ hence forward. Equations (\ref{Hf}) through (\ref{paramCond2}) must hold as well (with a minor modification to Eqs. (\ref{paramCond1})-(\ref{paramCond2}), cf Appendix), and lead to obtaining the set of boundary conditions collected in vector $\boldsymbol{Y}$, in the same form of Eqs. (\ref{equalirycontrF}) (case (a)) or (\ref{equalirycontrB}) (case (b)). 

A major challenge of multi-arc optimal control problems resides in the existence of additional multipoint corner conditions, which must be enforced at transition times \cite{pontani2021optimal}. However, for the problem at hand, previous research \cite{beolchi2024} proved that these conditions can be combined, to yield the implicit costate transformation,
\begin{equation}\label{CostateImplicit}
        \left( \frac{\partial \boldsymbol{\chi}_j}{\partial \boldsymbol{x}_{ini}^{(j + 1)}} \right)^{-T} \boldsymbol{\lambda}_{ini}^{{(j + 1)}} + \left( \frac{\partial \boldsymbol{\chi}_j}{\partial \boldsymbol{x}_{fin}^{(j)}} \right)^{-T} \boldsymbol{\lambda}_{fin}^{{(j)}} = \boldsymbol{0}.
\end{equation}
This relation allows finding $\boldsymbol{\lambda}_{ini}^{(j + 1)}$ from $\boldsymbol{\lambda}_{fin}^{(j)}$.  As a result, the set of unknown quantities reduces to the same set of a single-arc optimal control problem, i.e., the unknown quantities are the ones collected in $\boldsymbol{X}$ (cf. Eqs. (\ref{XvecF}) and (\ref{XvecB}), for cases (a) and (b), respectively).

\subsection{Method of solution}
The method of solution employs the indirect heuristic method, based on the joint use of the necessary conditions for an extremal, described in the preceding section and in the Appendix, and a heuristic technique, i.e., DE in this research. The same canonical units employed for transfers between LLO and Gateway are used.

 The numerical solution process utilizes a population of individuals, with fitness function reported in Eq. (\ref{Jtilde}). Based on the necessary conditions, the solution technique consists of the following steps: 

 \begin{enumerate}[itemsep=-0.1cm]
     \item Identify the known initial values of the state and costate variables and the minimal set of unknown values.
     \item For each individual $i$, with $i \leq N_p$, iterate the following sub-steps:
\begin{enumerate}
    \item select the unknown values for the quantities contained in $\boldsymbol{Y}$ (cf. Eqs. (\ref{equalirycontrF}) and (\ref{equalirycontrB})); 
    \item after selecting either $t_0$ or $t_f$, calculate the initial values of the state components;
    \item integrate numerically the state and costate equations while using Eqs. (\ref{anglesalfa})-(\ref{anglesbeta}) until Eq. (\ref{SOImoon}) is fulfilled, then go to sub-step (e); 
    \item if the preceding equality is never met, then set the fitness function $\tilde{J}$ to a large value and step to the next individual;
    \item apply four times the implicit costate transformation (\ref{CostateImplicit}), to get $\boldsymbol{\lambda}_{ini}^{(5)}$ from $\boldsymbol{\lambda}_{fin}^{(1)}$;
    \item apply four times the explicit state transformation (\ref{StateExplicitTransfTime}), to get $\boldsymbol{x}_{ini}^{(5)}$ from $\boldsymbol{x}_{fin}^{(1)}$;
    \item integrate numerically the state and costate equations while using Eqs. (\ref{anglesalfa})-(\ref{anglesbeta}) until $\tau = \tau_{fin}$;
    \item evaluate the Hamiltonian at the final time and the integral term that appears in the last component of $\boldsymbol{Y}$ (cf. Eqs. (\ref{equalirycontrF}) and (\ref{equalirycontrB}));
    \item if inequality (\ref{Hf}) is met, then go to sub-step (j), otherwise set $\tilde{J}$ to a large value and step to the next individual;
    \item evaluate the violations of the boundary conditions (cf. Eqs (\ref{equalirycontrF}) or (\ref{equalirycontrB})), then compute the fitness function $\tilde{J}$ (cf. Eq. (\ref{Jtilde})); 
    \item go to the next individual.
\end{enumerate}
     \item Once the fitness functions are obtained for all the individuals that compose the population, use the DE algorithm to update the position of individuals in the search space.
     \item Repeat all the steps until the fitness function of the best individual reaches a value lower than a threshold or until the maximum number of iterations is reached.
     
 \end{enumerate}

\indent DE stops when a prescribed number of iterations is reached or if the global best has not changed for several stalled iterations. In this work, the same settings adopted in \cite{beolchi2024} are used.

\subsection{Numerical results}
\indent The point-mass spacecraft is characterized by the propulsion parameters reported in Eq. (\ref{proppar}). The fitness function is $\Tilde{J}$ (cf. Eq. (\ref{Jtilde}), with the same weights used for Gateway-to-LLO transfers). The initial values of the adjoint variables are sought in the interval $[-1,1]$, due to their scalability, while the time of flight is sought in the interval [100,170] days, based on the analytical estimate described in \cite{beolchi2024}. In this approximate evaluation, the final lunar orbit is assumed to be polar and circular, with radius equal to the periselenium radius of Gateway. This approach reasonably provides an overestimate of the time of flight, because the final orbit is much closer to the Moon than the real path flown by Gateway. This rough estimate provides a time of flight of approximately 170 days, whereas the lower bound of the interval was set after an initial simulation session.

For case (a) (Gateway-to-LEO transfer), the spacecraft travels the NRHO flown by Gateway, and the initial time of the transfer (i.e., $t_0$) is selected along a single period of the departing orbit, in the interval [23 May 2025 at 22:35, 30 May 2025 at 05:38] UTC. The final orbit is a LEO with the following desired orbit elements: 
\begin{equation}\label{edid}
p_d = R_{E} + 463 \, \, {\rm{km}} \, \, \, \, \,  e_d = 0 \, \, \, \,  \,  i_d = 51.6^{\circ}
\end{equation}
where $R_{E} = 6378.136 \, \, \rm{km}$ is the mean equatorial radius of the Earth. The preceding values approximately correspond to the orbit traveled by the International Space Station.

\indent After running DE, the \textit{fminsearch} MATLAB routine is employed for refinement, leading to $\Tilde{J} = 5.012 \cdot 10^{-9}$. 
The optimal starting time is 29 May 2025 at 09:48:53 UTC and the time of flight equals $ \tau_{fin} = 153 \, \, {\rm d} \, \, 6 \, \, {\rm hrs} \, \, 51 \, \, {\rm min} \, \, 5 \, \, {\rm s} $, whereas the final mass ratio is $0.7835$. Table \ref{tabstate3} provides the spacecraft orbit elements at the initial and final time. The orbit elements at $t_0$ are the osculating elements of Gateway at the beginning of the transfer. \\

\begin{figure}[H]
    \includegraphics[width=1 \columnwidth]{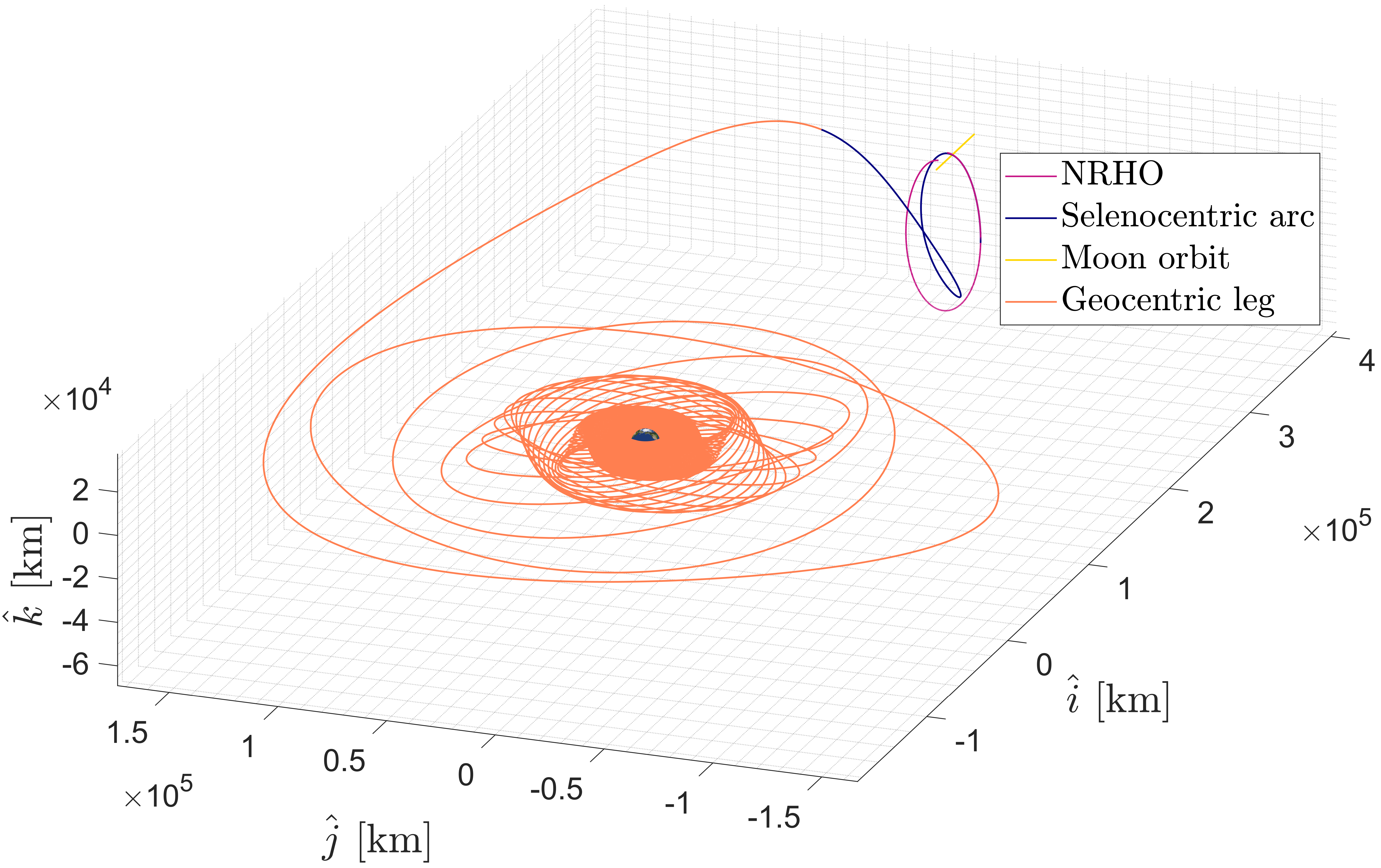}
    \caption{Gateway-to-LEO: minimum-time transfer in the synodic frame}
    \label{trajoptF2}
\end{figure}

\bgroup
\def\arraystretch{1}
\begin{table}[h] 
\centering
\begin{tabular}{ccc}
\toprule
\text{COE} & ${t_0}$ & ${t_f}$\\
\midrule
$a$ [km]&$4.418\cdot10^{4}$&$6.841\cdot10^{3}$\\ 
$e$&$0.927$ & $5.941\cdot10^{-6}$\\ 
$i\,[\degree]$  &$97.82 $ & $51.6$\\
$\Omega\,[\degree]$  &$ 20.56$ &$8.62 $ \\
$\omega\,[\degree]$  &$91.34 $ & $32.04 $ \\
$\theta_*\,[\degree]$  &$-155.58 $ & $86.37 $  \\
\bottomrule
\end{tabular}
\vspace*{1mm}
\caption{Gateway-to-LEO: initial and final orbit elements (referred to Moon and Earth, respectively)}
\label{tabstate3}
\end{table}
\egroup
\begin{figure}[H]
\centering
   \includegraphics[width=1\columnwidth]{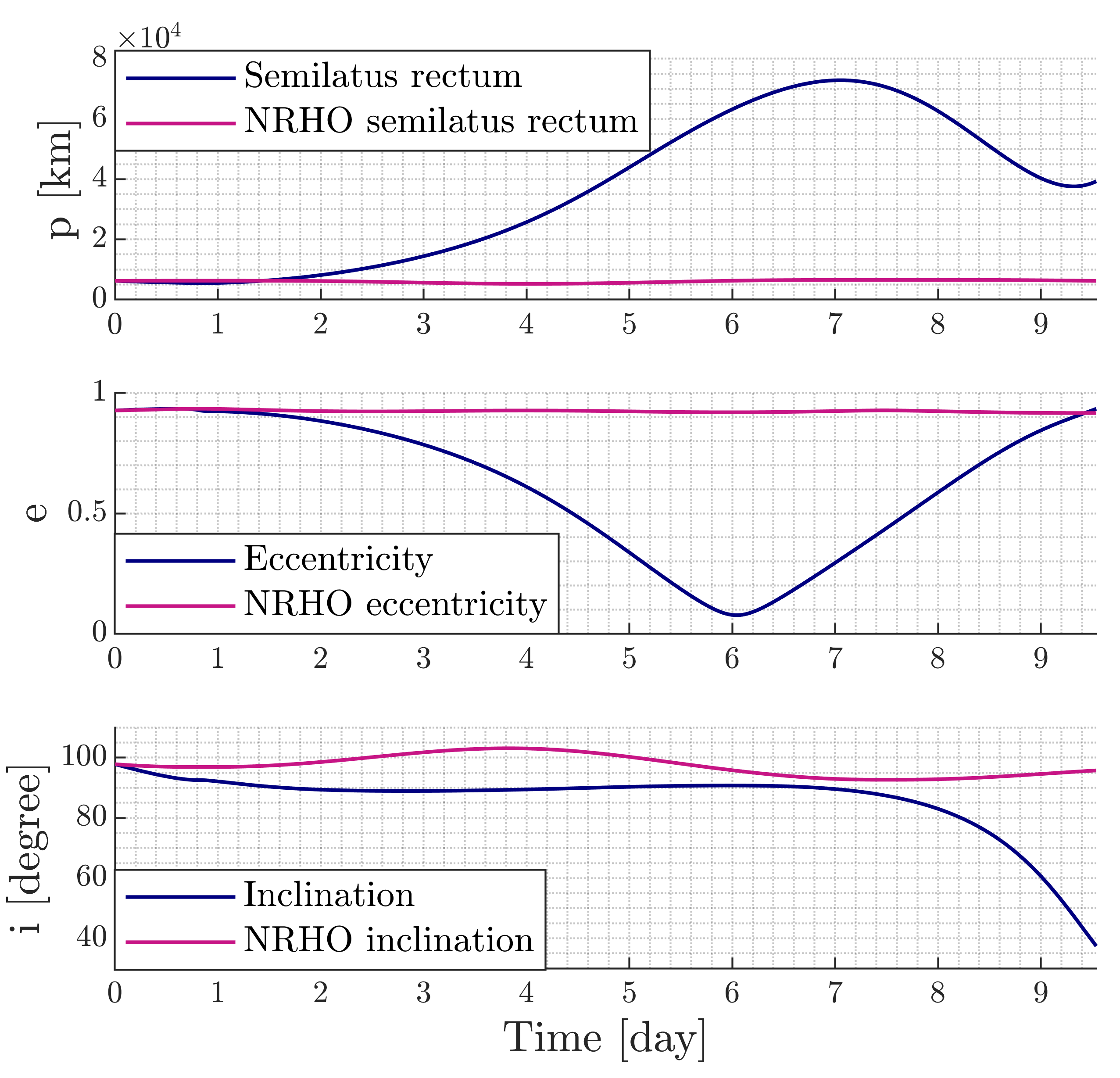}
   \caption{Gateway-to-LEO: time histories of the relevant orbit elements in the selenocentric arc}
    \label{pei_optF2}
\end{figure}
\begin{figure}[H]
\centering
     \includegraphics[width=1\columnwidth]{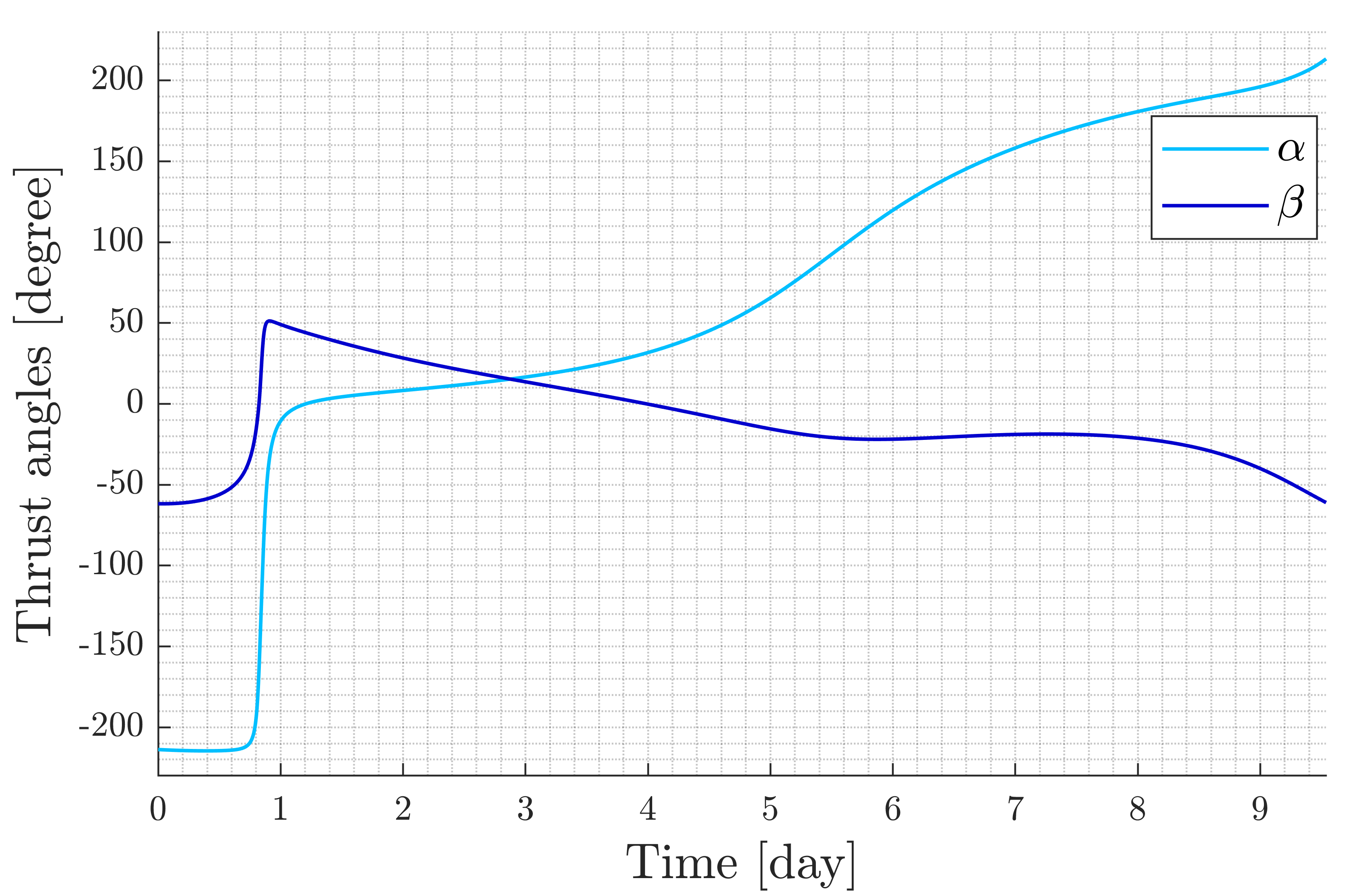}
    \caption{Gateway-to-LEO: time histories of thrust angles in the selenocentric arc}
    \label{alfabeta_optF2}
\end{figure}
From inspection of the preceding figures, it is apparent that in the selenocentric arc both the orbit elements and the thrust angles vary smoothly, without exhibiting oscillations. In contrast, in the geocentric arc, the thrust direction oscillates, while tending to be planar and opposed to the direction of motion. This braking effect corresponds to spiraling motion about Earth, with final injection into the desired circular LEO. This is also apparent from Fig. \ref{trajoptF2}.
\begin{figure}[H]
\centering
   \includegraphics[width=1\columnwidth]{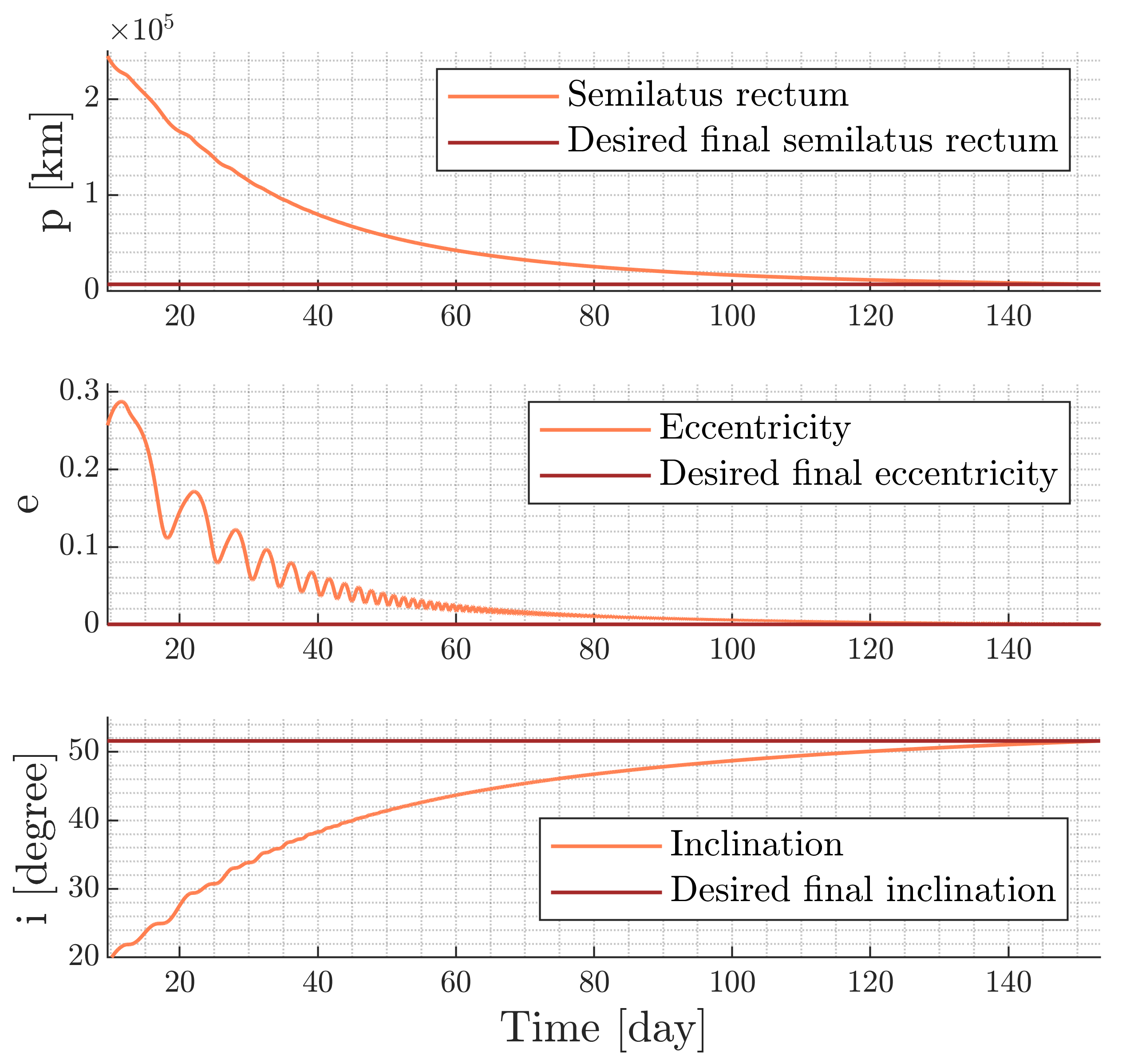}
   \caption{Gateway-to-LEO: time histories of the relevant orbit elements in the geocentric arc}
    \label{pei_optF3}
\end{figure}
\begin{figure}[H]
\centering
     \includegraphics[width=1\columnwidth]{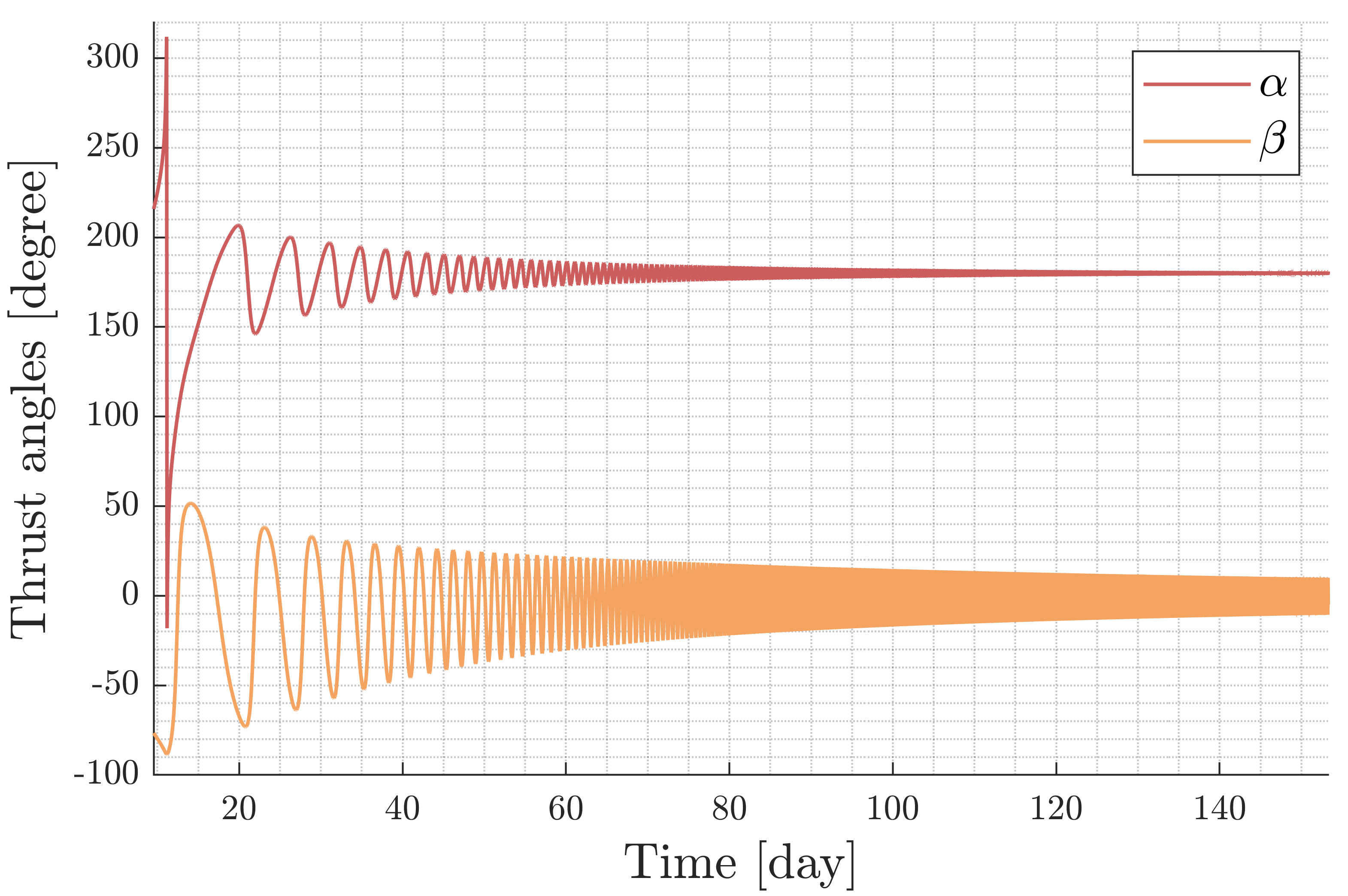}
    \caption{Gateway-to-LEO: time histories of thrust angles in the geocentric arc}
    \label{alfabeta_optF3}
\end{figure}
For case (b) (LEO-to-Gateway transfer), initially the spacecraft of interest travels a circular LEO with specified inclination and the final time of the transfer (i.e., $t_f$) is selected in the same time interval of case (a), referred to the period of the NRHO traveled by Gateway. More precisely, the orbit at $t_0$ is a LEO with orbit elements reported in Eq. (\ref{edid}).
\indent Figure \ref{pei_optF2} depicts the time histories of the lunar orbit elements $p$, $e$, and $i$, whereas Fig. \ref{alfabeta_optF2} portrays the optimal thrust angles in the selenocentric arc. Figures \ref{pei_optF3} and \ref{alfabeta_optF3} illustrate the respective variables in the geocentric arc. Finally, Fig.\ref{trajoptF2} depicts the minimum-time Gateway-to-LEO transfer, in the synodic frame. \\
\indent After running DE, the \textit{fminsearch} MATLAB routine is employed for refinement, leading to $\Tilde{J} = 1.112 \cdot 10^{-9}$. 
The optimal starting time is 05 January 2025 at 09:43:23 UTC and the time of flight equals $ \tau_{fin} = 144 \, \, {\rm d} \, \, 4 \, \, {\rm hrs} \, \, 50 \, \, {\rm min} \, \, 18 \, \, {\rm s} $, whereas the final mass ratio is $0.7963$. Table \ref{tabstate4} provides the spacecraft orbit elements at the initial and final time. The orbit elements of the spacecraft at $t_f$ are the osculating elements of Gateway at the end of the transfer. \\
\bgroup
\def\arraystretch{1.1}
\begin{table}[h] 
\centering
\begin{tabular}{ccc}
\toprule
\text{COE} & ${t_0}$ & ${t_f}$\\
\midrule
$a$ [km] &$6.841\cdot10^{3}$&$4.556\cdot10^{4}$\\ 
$e$ &$4.457\cdot10^{-6}$ & $0.929$\\ 
$i\,[\degree]$  &$51.6 $ & $97.39 $\\
$\Omega\,[\degree]$  &$9.48 $ &$ 21.40$ \\
$\omega\,[\degree]$  &$-183.92 $ & $91.10 $ \\
$\theta_*\,[\degree]$  &$-29.22 $ & $-151.20$  \\
\bottomrule
\end{tabular}
\vspace*{3mm}
\caption{LEO-to-Gateway: initial and final orbit elements (referred to Earth and Moon, respectively)}
\label{tabstate4}
\end{table}
\egroup

\begin{figure}[H]
   \includegraphics[width=1\columnwidth]{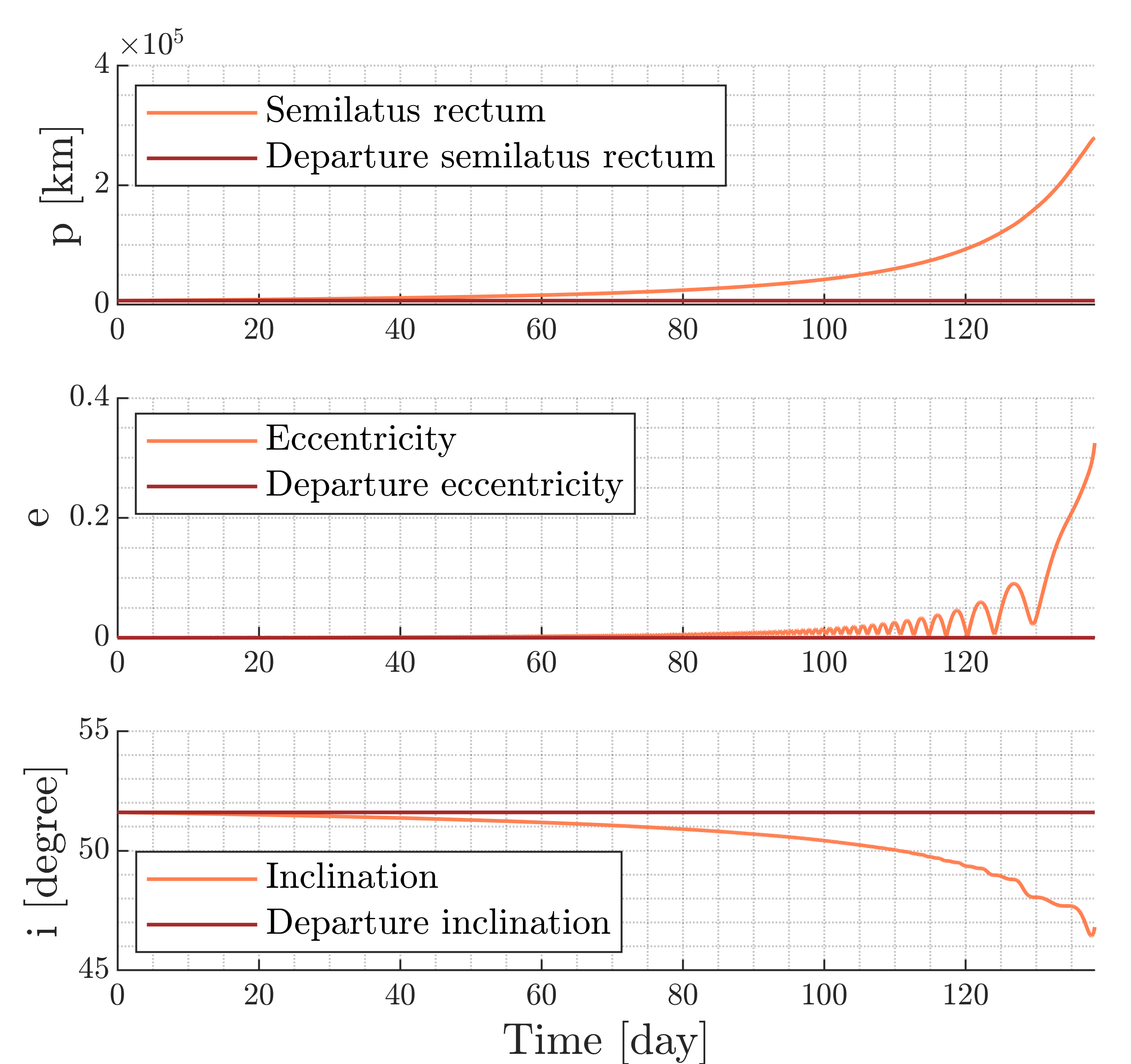}
   \caption{LEO-to-Gateway: time histories of the relevant orbit elementsin the geocentric arc}
    \label{pei_optF4}
\end{figure}
\begin{figure}[h]
\centering
     \includegraphics[width=1\columnwidth]{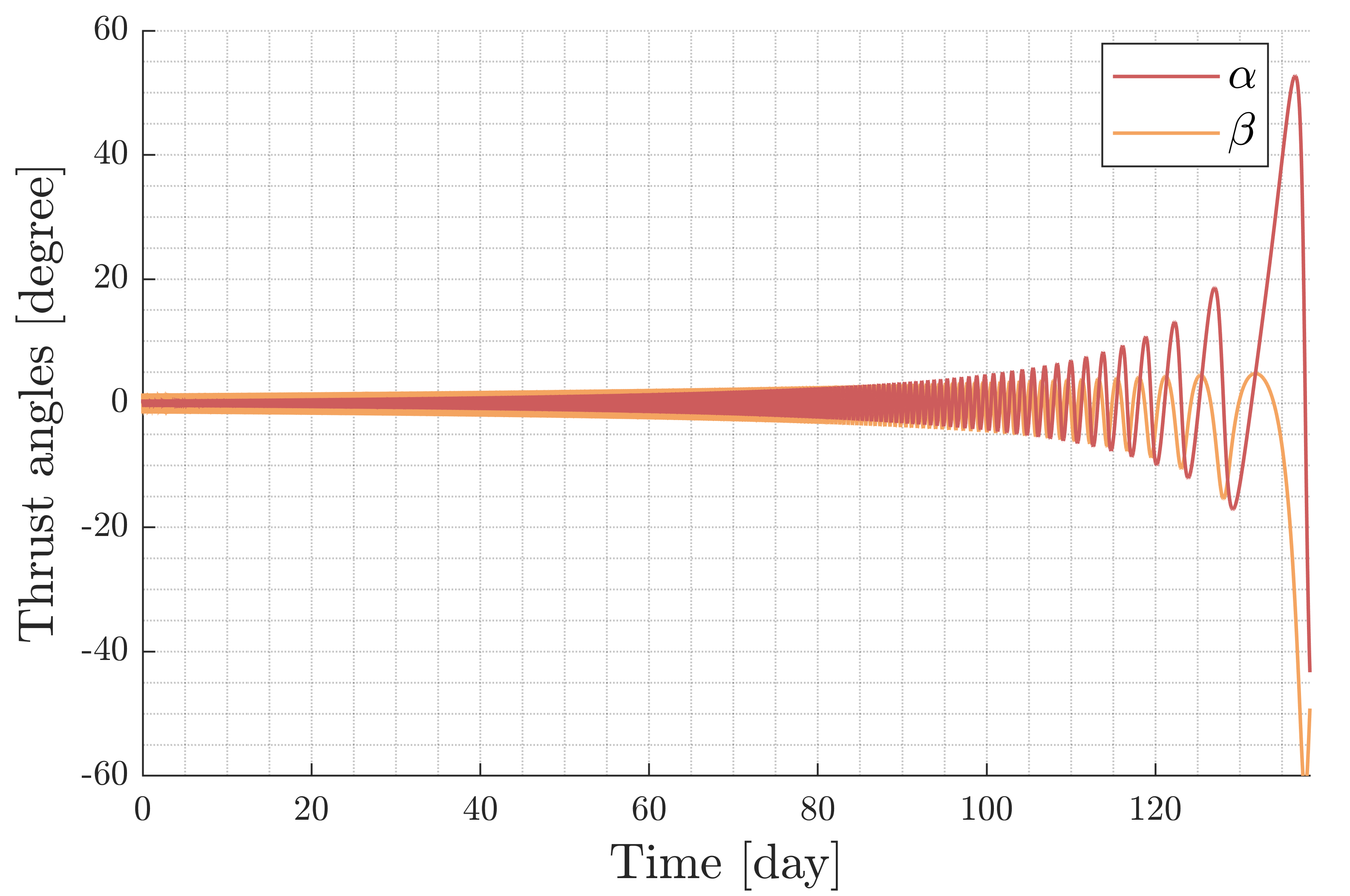}
    \caption{LEO-to-Gateway: time histories of thrust angles in the geocentric arc}
    \label{alfabeta_optF4}
\end{figure}
\begin{figure}[H]
\centering
   \includegraphics[width=1\columnwidth]{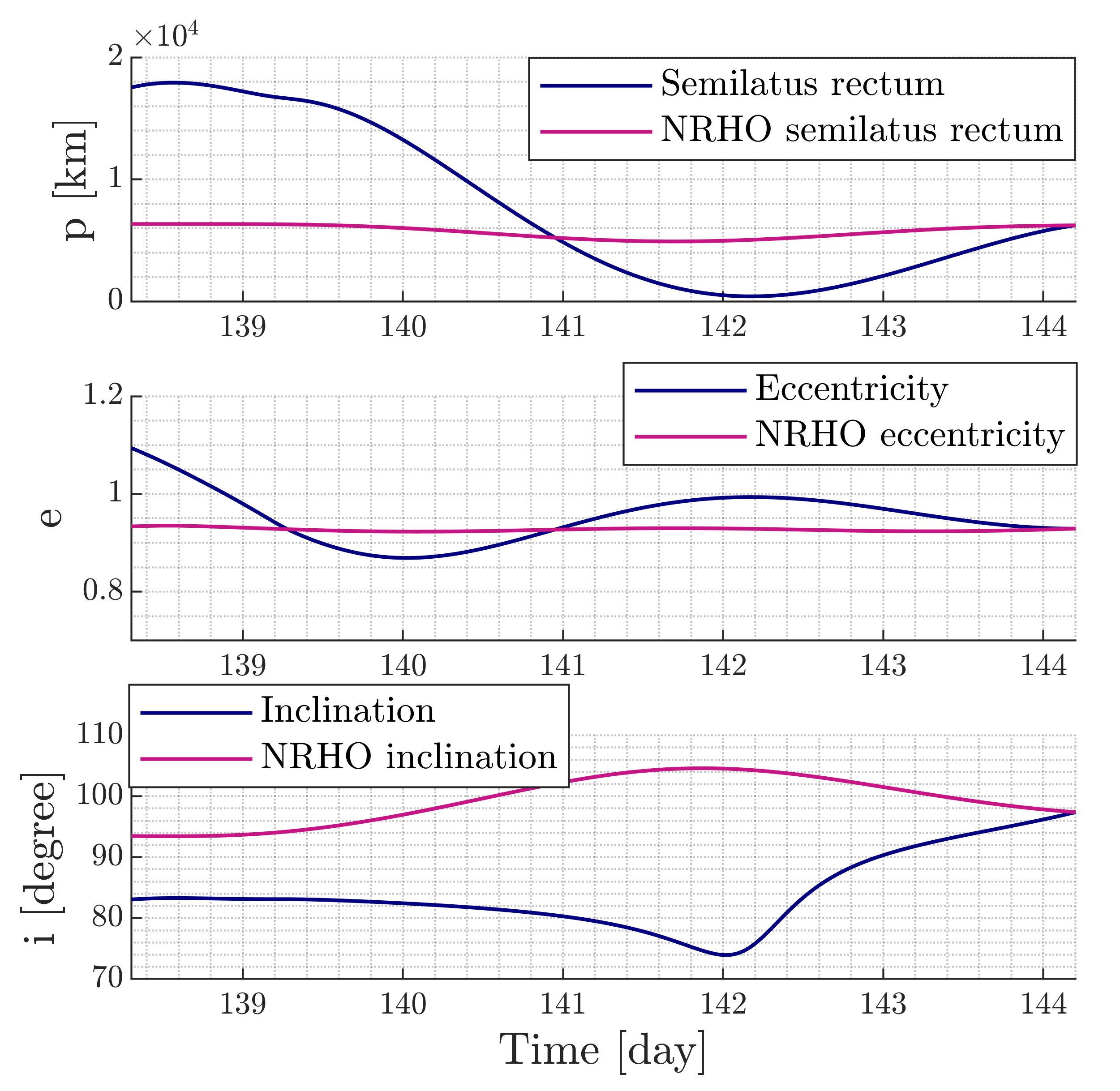}
   \caption{LEO-to-Gateway: time histories of the relevant orbit elements in the selenocentric arc}
    \label{pei_optF5}
\end{figure}
\begin{figure}[H]
\centering
     \includegraphics[width=1\columnwidth]{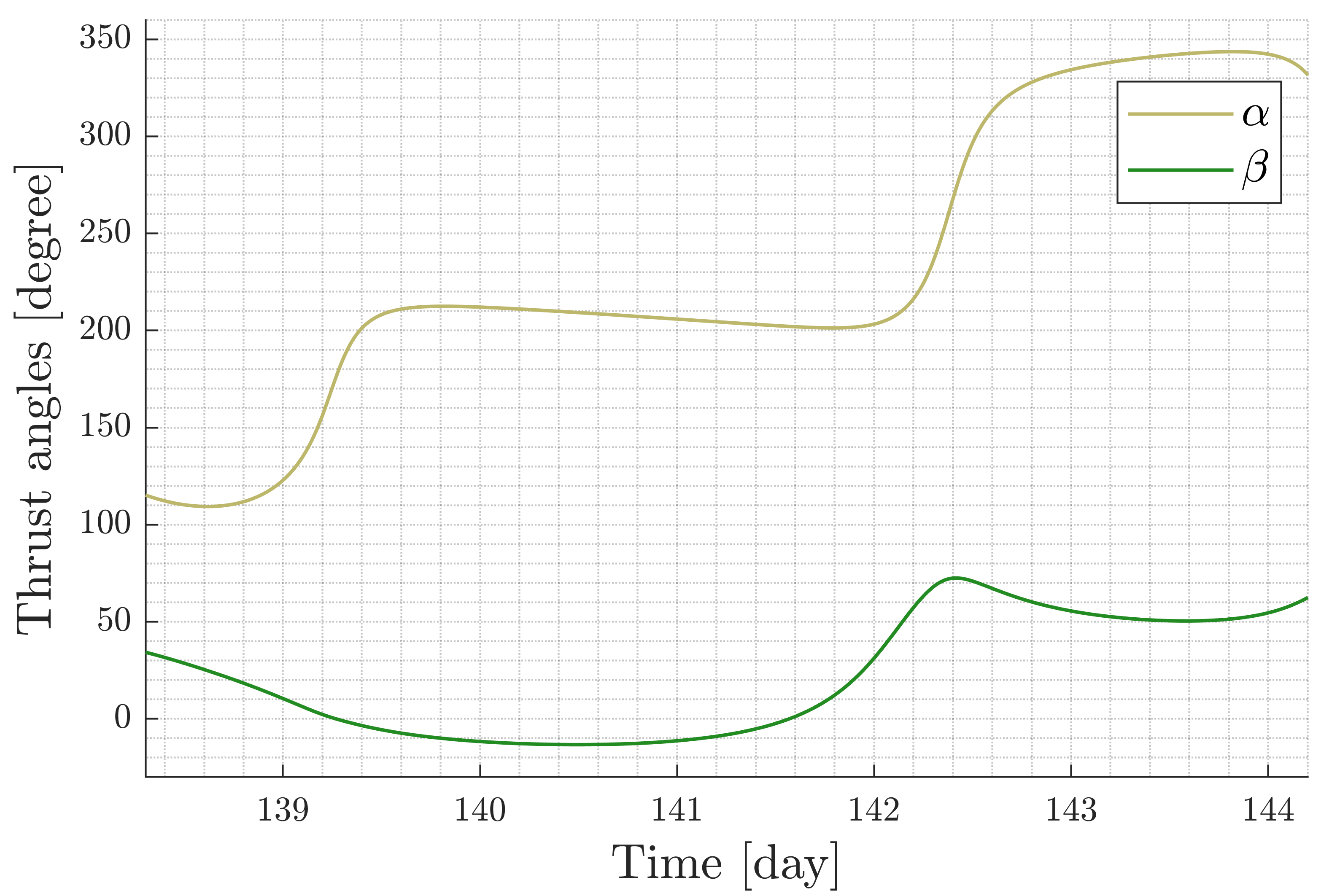}
    \caption{LEO-to-Gateway: time histories of thrust angles in the selenocentric arc}
    \label{alfabeta_optF5}
\end{figure}
\begin{figure}[H]
    \includegraphics[width=\columnwidth]{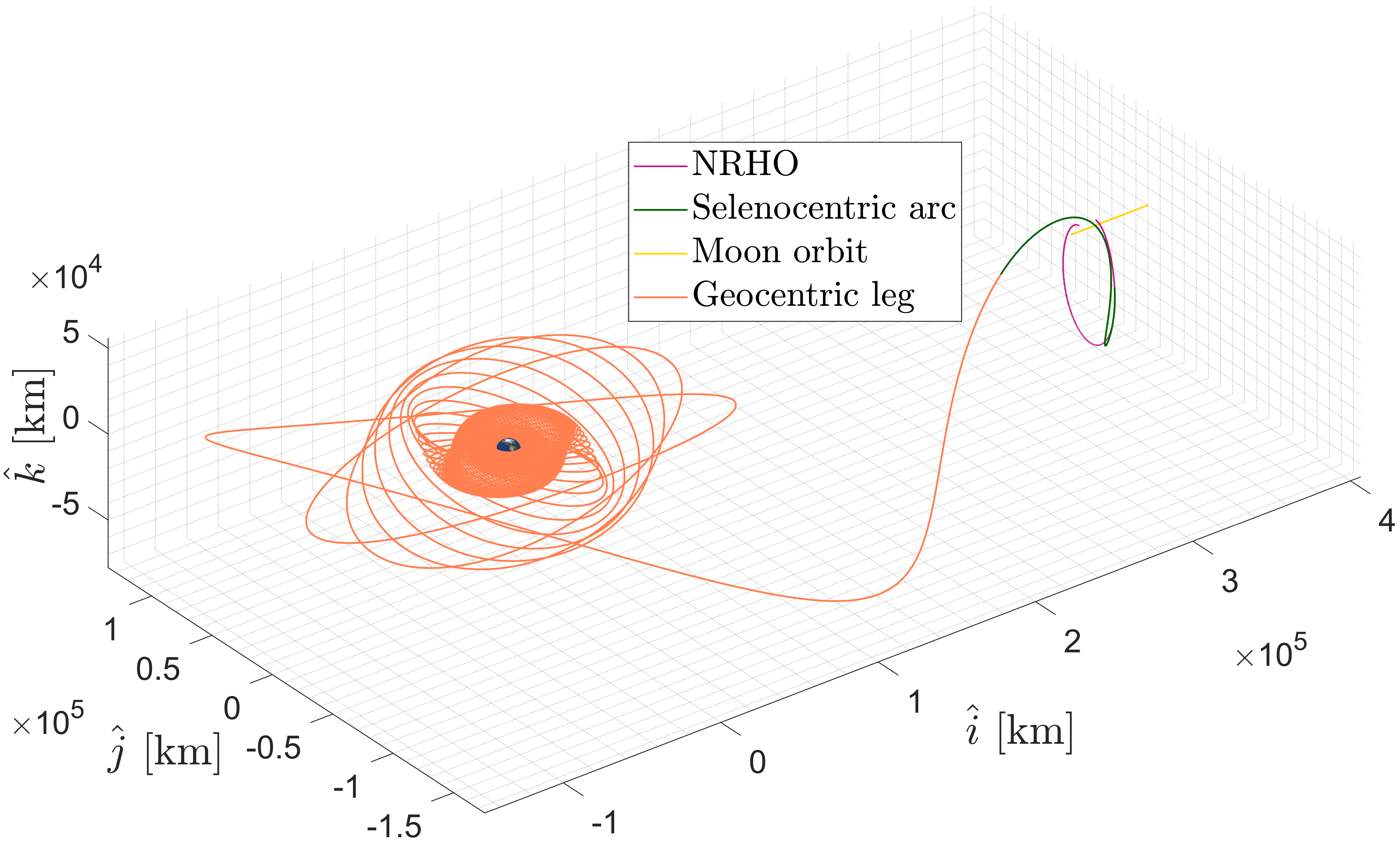}
    \caption{LEO-to-Gateway: minimum-time transfer in the synodic frame}
    \label{trajoptF3}
\end{figure}

\indent Figure \ref{pei_optF4} depicts the time histories of the terrestrial orbit elements $p$, $e$, and $i$, whereas Figs. \ref{alfabeta_optF4} portrays the optimal thrust angles in the geocentric arc. Figures \ref{pei_optF5} and \ref{alfabeta_optF5} illustrate the respective variables in the selenocentric arc. Finally, Fig.\ref{trajoptF3} depicts the minimum-time LEO-to-Gateway transfer, in the synodic frame.
From inspection of the preceding figures, it is apparent that spiraling motion characterizes the early phases in the geocentric arc, with a thrust direction that exhibits small oscillations about the direction of the instantaneous velocity, up to 100 days.
  In the selecnocentric arc, the thrust angles vary smoothly, as well as the osculating lunar orbit elements, and finally the spacecraft reaches Gateway, while matching all of its orbit elements (although only 3 of them are shown in Fig. \ref{pei_optF5}).


\section{CONCLUDING REMARKS}
\indent This study is devoted to determining minimum-time low-thrust orbit transfers able to connect both Earth and Moon to Gateway, which is placed in a lunar near-rectilinear halo orbit. 
Low-thrust orbit dynamics is propagated in a high-fidelity dynamical framework, with the use of planetary ephemeris and the inclusion of the simultaneous gravitational action of Sun, Earth, and Moon, along the entire transfer paths, in all cases. Modified equinoctial elements (MEE) are used because they avoid the hypersensitivity issues encountered with alternative state representations in indirect numerical solution approaches. Moreover, as a novel and effective approach with respect to preceding methodologies, backward propagation is proposed as a very convenient strategy for trajectories that approach and rendezvous with Gateway. This leads to developing an original, unified formulation of the underlying optimal control problem, capable of incorporating both forward and backward propagation, i.e., two-way transfer paths that originate from or arrive at Gateway. 
First, two-way transfers from Gateway to a specified lunar orbit are determined, using the indirect heuristic method, which employs the necessary conditions for optimality and a heuristic algorithm. Second, orbit transfers that connect Earth and Gateway are addressed. These trajectories are affected by two major attracting bodies, and the related optimal control problem is formulated and solved in a multibody dynamical framework. Unlike most previous contributions focused on indirect approaches applied to optimal low-thrust trajectories, this research extends the indirect formulation of minimum-time orbit transfers to a multibody mission scenario, with terminal orbits around two distinct primaries, i.e., Earth and Moon. Three different representations for the spacecraft state are employed, i.e. (i) MEE relative to Earth, (ii) Cartesian coordinates, and (iii) MEE relative to Moon, and the multi-arc trajectory optimization problem at hand includes two legs: (a) geocentric leg and (b) selenocentric leg. While MEE are employed for orbit propagation, Cartesian coordinates are useful as intermediate variables, to step from Earth MEE to lunar MEE (and viceversa). Multi-arc optimal control problems are associated with several additional corner conditions, at the transition between the two legs. This study derives a closed-form solution to these corner conditions, using implicit costate transformation. As a result, the parameter set for an indirect algorithm retains the reduced size of the typical set associated with a single-arc optimization problem and the indirect heuristic technique is applied again. Moreover, the use of backward propagation for transfer paths that travel toward Gateway allows considerable simplification to the numerical solution technique, which does not require any stratified objective function, unlike what occurred in previous research. The numerical results unequivocally prove that the methodology developed in this study is effective for determining minimum-time low-thrust orbit transfers connecting Earth, Moon, and Gateway, in a high-fidelity ephemeris-based dynamical framework.

\section*{{APPENDIX A}}

\indent This appendix is concerned with the formulation and necessary conditions for an extremal of forward and backward minimum-time low-thrust orbit transfers. This is especially useful for orbit transfers that involve departure from or arrival at Gateway. In fact, both departure from or rendezvous with Gateway require exact matching of its position and velocity. Therefore, it is convenient to introduce the auxiliary independent variable $\tau$, which plays the role of either forward time (for transfers that originate from Gateway) or backward time (for transfers that rendezvous with Gateway). In both cases, the state equations have the form 
\begin{equation}\label{stateEqTau}
\dot{\boldsymbol{x}}=\boldsymbol{f}\left(\boldsymbol{x},\boldsymbol{u},\boldsymbol{p},\tau\right)       
    \end{equation}
where
\begin{equation}\label{tauF}
\tau=t-t_0\quad\textrm{for forward transfers}       
    \end{equation}
or
\begin{equation}\label{tauB}
\tau=t_f-t\quad\textrm{for backward transfers}       
    \end{equation}
In both cases, $\tau$ is monotonically increasing and constrained to $[0,t_f-t_0]$. The two epochs $t_0$ and $t_f$ are collected as components of the parameter vector $\boldsymbol{p}:=\left[p_1\,\,p_2\right]^T=\left[t_0\,\,t_f\right]^T$, and the following equality constraint is introduced:
    \begin{equation}\label{taufinAndp}
\tau_{fin}+p_1-p_2=0        
    \end{equation}
    Hence, the vector of boundary conditions assumes the form
\begin{equation}\label{Psi}
    \boldsymbol{\Psi} \left( \boldsymbol{x}_{i}, \boldsymbol{x}_{f}, \tau_{i}, \tau_{fin}, \boldsymbol{p} \right) = \begin{bmatrix}
        \boldsymbol{\Psi}_i \left( \boldsymbol{x}_{i}, \tau_{i}, \boldsymbol{p}\right)\\ \boldsymbol{\Psi}_f \left( \boldsymbol{x}_{f}, \tau_{fin}, \boldsymbol{p} \right)\\
        \tau_{fin}+p_1-p_2
    \end{bmatrix} = \boldsymbol 0. 
\end{equation}
The objective function to minimize is the time of flight, i.e.,
\begin{equation}\label{objFun}
    J=k_J \tau_{fin}
\end{equation}
where $k_J$ denotes an arbitrary positive constant.

\indent The necessary conditions for an extremal derive from the calculus of variations, after introducing the Hamiltonian function $H$ and the function of boundary conditions $\Phi$ \cite{hull2013optimal},
\begin{equation}\label{Ham}
    H=\boldsymbol{\lambda}^T\boldsymbol{f}
\end{equation}
\begin{equation}\label{bigPhi}
    \Phi=k_J \tau_{fin}+\boldsymbol{\upsilon}_i^T\boldsymbol{\Psi_i}+\boldsymbol{\upsilon}_f^T\boldsymbol{\Psi_f}+\xi\left(\tau_{fin}+p_1-p_2\right)
\end{equation}
where $\boldsymbol{\lambda}$ represents the adjoint vector conjugate to the equations of motion, whereas $\boldsymbol{\upsilon}_i$, $\boldsymbol{\upsilon}_f$, and $\xi$ denote the adjoint variables associated with the boundary conditions. The necessary conditions for an extremal include the well-established adjoint equations and the Pontryagin minimum principle (omitted for the sake of conciseness), together with the boundary conditions on the adjoint vector, 
\begin{equation}\label{lambda0}
\boldsymbol{\lambda}_i=-\left(\frac{\partial\boldsymbol{\Psi}_i}{\partial \boldsymbol{x}_i}\right)^T \boldsymbol{\upsilon}_i
\end{equation}
\begin{equation}\label{lambdaf}
\boldsymbol{\lambda}_f=\left(\frac{\partial\boldsymbol{\Psi}_f}{\partial \boldsymbol{x}_f}\right)^T \boldsymbol{\upsilon}_f
\end{equation}
the transversality condition on the final Hamiltonian
\begin{equation}\label{Hf}
H_f+k_J-\xi=0\quad\textrm{i.e.,}\quad H_f-\xi=-k_J<0
\end{equation}
and the parameter condition \cite{hull2013optimal}
\begin{equation}\label{paramCond}
\int_{\tau_i}^{\tau_{fin}}\left(\frac{\partial H}{\partial\boldsymbol{p}}\right)^T d\tau+
\left(\frac{\partial \Phi }{\partial \boldsymbol{p}} \right)^T=\boldsymbol{0},\,\,\textrm{with}\,\,\tau_i=0,
\end{equation}
leading to
\begin{equation}\label{paramCond1}
\int_{0}^{\tau_{fin}}\left(\frac{\partial H}{\partial p_1}\right) d\tau+
\left(\frac{\partial\boldsymbol{\Psi_i}}{\partial p_1}\right)^T\boldsymbol{\upsilon}_i+\left(\frac{\partial\boldsymbol{\Psi_f}}{\partial p_1}\right)^T\boldsymbol{\upsilon}_f+\xi=0
\end{equation}
\begin{equation}\label{paramCond2}
\int_{0}^{\tau_{fin}}\left(\frac{\partial H}{\partial p_2}\right) d\tau+
\left(\frac{\partial\boldsymbol{\Psi_i}}{\partial p_2}\right)^T\boldsymbol{\upsilon}_i+\left(\frac{\partial\boldsymbol{\Psi_f}}{\partial p_2}\right)^T\boldsymbol{\upsilon}_f-\xi=0
\end{equation}
The last two relations represent a linear system of two equations in the unknowns $\boldsymbol{\upsilon}_i$, $\boldsymbol{\upsilon}_f$, and $\xi$. The Hamiltonian depends on $\boldsymbol{p}$ through the perturbing action of third bodies, related to the epoch (i.e., $t_0$ and $t_f$, cf. also Appendix B). It is worth remarking that the final condition on the Hamiltonian (\ref{Hf}) is equivalent to an inequality constraint, due to arbitrariness of $k_J$. In contrast, no condition exists on the initial Hamiltonian, because $\tau_i$ is specified and equal to 0.

It is relatively straightforward to conclude that the necessary conditions are well posed and, in conjunction with the state equations (\ref{stateEqTau}) and the boundary conditions (\ref{Psi}), correspond to the set of unknown variables, namely the state and control vectors $\boldsymbol{x}$ and $\boldsymbol{u}$, the parameter vector $\boldsymbol{p}$, the time of flight $\tau_{fin}$, the adjoint variables $\boldsymbol{\lambda}$, $\boldsymbol{\upsilon}_i$, $\boldsymbol{\upsilon}_f$, and $\xi$. 

The preceding developments refer to a single-arc optimal control problem. Multi-arc optimal control problems require introducing additional corner conditions, leading to Eq. (\ref{CostateImplicit}), whereas the relations (\ref{Psi}), (\ref{objFun}), (\ref{bigPhi}) and (\ref{lambda0})-(\ref{Hf}) remain unaltered. Instead, Eqs. (\ref{paramCond})-(\ref{paramCond2}) must be modified, by replacing the integral term with a summation of integral terms. In particular, Eq. (\ref{paramCond}) must include a summation of integral terms with integrand $\partial H^{(j)}/\partial\boldsymbol{p}$ (cf. Eq. (\ref{Jbarmultiarc})). Similar modifications apply to Eqs. (\ref{paramCond1})-(\ref{paramCond2}).

\section*{{APPENDIX B}}

\indent This appendix is devoted to analyzing the parameter conditions that appear as the last rows of vector $\boldsymbol{Y}$ in Eqs. (\ref{equalirycontrF}) and (\ref{equalirycontrB}).

As a first step, the Hamiltonian can be recognized to depend explicitly on the actual time $t$, due to the perturbing actions of third bodies. In fact, depending on the mission scenario, Sun and either Earth or Moon play the role of third bodies, and their instantaneous position (needed in Eq. (\ref{a3B})) is obtained through interpolation of ephemeris-based data. This circumstance implies that the Hamiltonian $H$ depends explicitly on the actual time $t$.

For the sake of simplicity, a single-arc problem (suitable for transfers between Gateway and LLO) is considered in the following.

If forward propagation is employed, then $t=t_0+\tau$, therefore $H$ can be rewritten as an explicit function of $(t_0+\tau)$. This implies that 
\begin{equation}\label{partialH_p1_F}
\frac{\partial H}{\partial p_1}=\frac{\partial H}{\partial t_0}=\frac{\partial H}{\partial \tau}
\end{equation}
However, because the state and costate vectors obey the Hamilton equations, 
\begin{equation}\label{partialH_p1_F_bis}
\frac{\partial H}{\partial \tau}=\frac{d H}{d \tau}
\end{equation}
As a result, the relation that involves $\partial H / \partial p_1$ in Eq. (\ref{equalirycontrF}) becomes
\begin{equation}\label{newparamCondForw}
H_f-H_0+\left(\frac{\partial\boldsymbol{\Psi_i}}{\partial p_1}\right)^T\boldsymbol{\upsilon}_i+\left(\frac{\partial\boldsymbol{\Psi_f}}{\partial p_1}\right)^T\boldsymbol{\upsilon}_f+\xi=0
\end{equation}
Moreover, because $\boldsymbol{\Psi_i}$ and $\boldsymbol{\Psi_f}$ are independent of $p_2$, Eq. (\ref{paramCond2}) yields $\xi=0$. Furthermore,  $\boldsymbol{\Psi_f}$ is also independent of $p_1$. Therefore, Eq. (\ref{newparamCondForw}) simplifies to 
\begin{equation}\label{newparamCondForw_simple}
H_f-H_0+\left(\frac{\partial\boldsymbol{\Psi_i}}{\partial p_1}\right)^T\boldsymbol{\upsilon}_i=0
\end{equation}
Let vector $\boldsymbol{x}_G$ contain the MEE of Gateway at $t_0$ (cf. Eq. (\ref{Psi_0_a})). Using Eqs. (\ref{Psi_0_a}) and (\ref{lambdainia}), Eq. (\ref{newparamCondForw_simple}) becomes
\begin{equation}\label{newparamCondForw_simple2}
H_f-H_0+\boldsymbol{\lambda}_i^T \frac{\partial \boldsymbol{x}_G}{\partial t_0}=0
\end{equation}
leading to 
\begin{equation}\label{newparamCondForw_simple3}
H_0<\boldsymbol{\lambda}_i^T \frac{\partial \boldsymbol{x}_G}{\partial t_0}
\end{equation}
because $H_f<0$ from Eq. (\ref{Hf}), in which $\xi=0$. The preceding relation can be proven to be always satisfied at a minimizing solution. In fact, using Eqs. (\ref{Hsplit})--(\ref{anglesbeta}), at time $t_0$ the Hamiltonian function is given by
\begin{equation}\label{initH}
H_0=\boldsymbol{\lambda}_i^T \boldsymbol{f}=H_x-\sqrt{H_r^2+H_{\theta}^2+H_h^2}
\end{equation}
As the Gateway motion is subject to natural orbit dynamics, because $d \tau = d t_0$, the right-hand side in Eq. (\ref{newparamCondForw_simple3}) is recognized to equal the Hamiltonian function at $t_0$ in the absence of propulsion, i.e.,
\begin{equation}\label{inStateGateway}
\boldsymbol{\lambda}_i^T \frac{\partial \boldsymbol{x}_G}{\partial t_0}=\boldsymbol{\lambda}_i^T \frac{\partial \boldsymbol{x}_G}{\partial \tau}\Big|_{\tau=0}=H_x
\end{equation}
Combination of Eqs. (\ref{initH}) and (\ref{inStateGateway}) proves that the condition (\ref{newparamCondForw_simple3}) holds as a consequence of the Pontryagin minimum principle. 

If backward propagation is employed, then $t=t_f-\tau$, therefore $H$ can be rewritten as an explicit function of $(t_f-\tau)$. This implies that 
\begin{equation}\label{partialH_p1_F}
\frac{\partial H}{\partial p_2}=\frac{\partial H}{\partial t_f}=-\frac{\partial H}{\partial \tau}
\end{equation}
However, because the state and costate vectors obey the Hamilton equations, Eq. (\ref{partialH_p1_F_bis}) holds. As a result, the relation that involves $\partial H / \partial p_2$ in Eq. (\ref{equalirycontrB}) becomes
\begin{equation}\label{newparamCondBack}
H_0-H_f+\left(\frac{\partial\boldsymbol{\Psi_i}}{\partial p_2}\right)^T\boldsymbol{\upsilon}_i+\left(\frac{\partial\boldsymbol{\Psi_f}}{\partial p_2}\right)^T\boldsymbol{\upsilon}_f+\xi=0
\end{equation}
Moreover, because $\boldsymbol{\Psi_i}$ and $\boldsymbol{\Psi_f}$ are independent of $p_1$, Eq. (\ref{paramCond1}) yields $\xi=0$. Furthermore,  $\boldsymbol{\Psi_f}$ is also independent of $p_2$. Therefore, Eq. (\ref{newparamCondBack}) simplifies to 
\begin{equation}\label{newparamCondBack_simple}
H_0-H_f+\left(\frac{\partial\boldsymbol{\Psi_i}}{\partial p_2}\right)^T\boldsymbol{\upsilon}_i=0
\end{equation}
Let vector $\boldsymbol{x}_G$ contain the MEE of Gateway at $t_f$ (cf. Eq. (\ref{Psi_0_b})). Using Eqs. (\ref{Psi_0_b}) and (\ref{lambdainia}), Eq. (\ref{newparamCondBack_simple}) becomes
\begin{equation}\label{newparamCondback_simple2}
H_0-H_f+\boldsymbol{\lambda}_i^T \frac{\partial \boldsymbol{x}_G}{\partial t_f}=0
\end{equation}
leading to 
\begin{equation}\label{newparamCondBack_simple3}
H_0<-\boldsymbol{\lambda}_i^T \frac{\partial \boldsymbol{x}_G}{\partial t_f}
\end{equation}
because $H_f<0$ from Eq. (\ref{Hf}), in which $\xi=0$. The preceding relation can be proven to be always satisfied at a minimizing solution. In fact, using Eqs. (\ref{Hsplit})--(\ref{anglesbeta}), at time $t_f$, associated with $\tau=0$, the Hamiltonian function is given by
\begin{equation}\label{initH2}
H_0=\boldsymbol{\lambda}_i^T \boldsymbol{f}=H_x-\sqrt{H_r^2+H_{\theta}^2+H_h^2}
\end{equation}
As the Gateway motion is subject to natural orbit dynamics, because $d \tau = -d t_f$, the right-hand side in Eq. (\ref{newparamCondBack_simple3}) is closely related to the Hamiltonian function at $t_f$ in the absence of propulsion, i.e.,
\begin{equation}\label{inStateGateway2}
-\boldsymbol{\lambda}_i^T \frac{\partial \boldsymbol{x}_G}{\partial t_f} = \boldsymbol{\lambda}_i^T \frac{\partial \boldsymbol{x}_G}{\partial \tau}\Big|_{\tau=0}=H_x
\end{equation}
Combination of Eqs. (\ref{initH2}) and (\ref{inStateGateway2}) proves that the condition (\ref{newparamCondBack_simple3}) holds as a consequence of the Pontryagin minimum principle. 

The preceding developments refer to a single-arc optimal control problem. Extension to the multi-arc problem addressed in this research is relatively straightforward, due to continuity of the Hamiltonian function at junction points.

\section*{ACKNOWLEDGMENTS}

\indent C. Pozzi, A. Beolchi, and E. Fantino acknowledge Khalifa University of Science and Technology's internal grant CIRA-2021-65/8474000413. A. Beolchi and E. Fantino have received support 
from project 8434000368 (Khalifa University of Science and Technology's 6U Cubesat mission). E. Fantino acknowledges grant ELLIPSE/8434000533 (Abu Dhabi's Technology Innovation Institute) and projects PID2020-112576GB-C21 and PID2021-123968NB-I00 of the Spanish Ministry of Science and Innovation.

\bibliographystyle{ieeetr}
\bibliography{main.bib}

\end{document}